\documentclass[aos,citesort,seceqn,dvips]{arximspdf}

\doi{10.1214/10-AOS836}
\volume{39}
\issue{1}
\pubyear{2011}
\firstpage{201}
\lastpage{231}

\makeatletter
\newcommand{\R}{\mathbb{R}}

\newtheorem{theorem}{Theorem}

\newproclaim{condition}{Condition}

\newtheorem{corollary}{Corollary}

\newproclaim{definition}{Definition}

\newtheorem{lemma}{Lemma}

\newtheorem{proposition}{Proposition}

\makeatother

\begin{document}
\begin{frontmatter}

\title{Global uniform risk bounds for wavelet deconvolution estimators}

\runtitle{Wavelet density deconvolution}

\begin{aug}
\author[A]{\fnms{Karim} \snm{Lounici}\ead[label=e1]{k.lounici@statslab.cam.ac.uk}}
\and
\author[A]{\fnms{Richard} \snm{Nickl}\ead[label=e2]{r.nickl@statslab.cam.ac.uk}\corref{}}

\runauthor{K. Lounici and R. Nickl}

\affiliation{University of Cambridge}
\address[A]{Statistical Laboratory\\
Department of Pure Mathematics\\
\quad and Mathematical Statistics\\
University of Cambridge\\
Cambridge, CB3 0WB\\
United Kingdom\\
\printead{e1}\\
\phantom{E-mail:\ }\printead*{e2}} 
\end{aug}

\received{\smonth{1} \syear{2010}}
\revised{\smonth{5} \syear{2010}}

\begin{abstract}
We consider the statistical deconvolution problem where one observes
$n$ replications from the model $Y=X+\epsilon$, where $X$ is the
unobserved random signal of interest and  $\epsilon$ is an independent
random error with distribution~$\varphi$. Under weak assumptions on the
decay of the Fourier transform of $\varphi,$ we derive upper bounds for
the finite-sample sup-norm risk of wavelet deconvolution density
estimators $f_n$ for the density $f$ of $X$, where $f\dvtx \mathbb R
\to \mathbb R$ is assumed to be bounded. We then derive lower bounds
for the minimax sup-norm risk over Besov balls in this estimation
problem and show that wavelet deconvolution density estimators attain
these bounds. We further show that linear estimators adapt to the
unknown smoothness of $f$ if the Fourier transform of $\varphi$ decays
exponentially and that a corresponding result holds true for the hard
thresholding wavelet estimator if $\varphi$ decays polynomially. We
also analyze the case where $f$ is a ``supersmooth''/analytic density.
We finally show how our results and recent techniques from Rademacher
processes can be applied to construct global confidence bands for the
density $f$.
\end{abstract}


\end{frontmatter}

\section{Introduction}

Consider the statistical deconvolution model
\begin{equation}\label{deconvolution-model}
Y = X + \epsilon,
\end{equation}
where $X$ is a real-valued random variable with unknown probability
density $f :  \mathbb{R}\rightarrow \mathbb{R}^{+}$ and $\epsilon$
is an error term independent of $X$ that is distributed according to
the probability measure $\varphi$ on $\mathbb R$. The law $P$ of $Y$
equals the convolution $f\ast\varphi$ and we denote its density by $g$.
Let $Y_{1},\ldots,Y_{n}$ be i.i.d.~replications of $Y$ in the model
(\ref{deconvolution-model}) and denote by $P_n$ the associated
empirical measure. The \textit{deconvolution problem} is about
recovering the unknown density $f$ from the noisy observations
$(Y_{1},\ldots,Y_{n})$. It has been extensively studied: we refer to
Carroll and Hall \cite{CarrollHall1988}, Stefanski \cite{Stefanski1990},
Stefanski and Carroll \cite{SC90}, Fan \cite{F91,fan1993}, Diggle and Hall \cite{DH93},
Goldenshluger \cite{Gold1999}, Pensky and Vidakovic \cite{PV99}, Delaigle and Gijbels
\cite{DG04}, Hesse and Meister \cite{HM04}, Johnstone et~al. \cite{JKPR04}, Johnstone and Raimondo \cite{JR04}, Bissantz et~al. \cite{BDHM07}, Bissantz and Holzmann \cite{BH}, Meister \cite{M08},
Butucea and Tsybakov \cite{BT08I,BT08II} and Pensky and Sapatinas \cite{PS09},
also to the monograph Meister \cite{M09}, as well as to Cavalier \cite{cavalier2008} for
a survey of the literature on general inverse problems in statistics,
of which deconvolution is a special case.

One key lesson from the aforementioned literature is that a lower bound
on the regularity of the signal $\epsilon$ is necessary to be able to
estimate $f$ with reasonable accuracy. This lower bound is often
quantified by a lower bound on the decay of the Fourier transform
$F[\varphi]$ of $\varphi$ and Fourier inversion techniques are applied
to construct estimators for $f$.

Most of the literature on this problem (with some notable exceptions,
to be discussed below) deals with the $L^2$-theory, that is, involves
the loss function $d^2(\hat f,f)=\int(\hat f-f)^2$ and is often
restricted to the case of periodic and hence compactly supported $f$.
These restrictions are theoretically convenient, in particular since
Fourier analysis-based methods can be used without too much difficulty,
using the Parseval--Plancherel isometry. However, a sound understanding
of the local behavior of deconvolution estimators seems to be of
significant statistical importance. In particular a theory that could
deal with sup-norm loss $d(\hat f,f)=\sup_{x \in \mathbb R}|\hat
f(x)-f(x)|$ could be used in the construction of confidence bands for
the object $f$ of statistical interest. A fortiori it is not at all
clear whether the intuitions from $L^2$-theory carry over to pointwise
and uniform loss functions in generality, bearing in mind that
$L^2$-convergence properties of Fourier series can give a completely
inadequate picture of their pointwise or uniform behaviour.

In the present article, we use methods from empirical process theory to
derive finite-sample sup-norm risk bounds for deconvolution density
estimators based on Fourier inversion with Meyer (or similar
band-limited) wavelets. These estimators were studied in Pensky and
Vidakovic \cite{PV99} and Johnstone et al. \cite{JKPR04}, and have since been
successfully used in inverse problems. Our results hold under minimal
assumptions on the density $f$ and the distribution $\varphi$: we
require $f$ to be bounded, which is unavoidable if one considers
sup-norm loss, and we assume that the Fourier transform of $\varphi$ is
nonzero on the intervals of support of the Meyer wavelet, which is
necessary to define any estimator based on Fourier inversion and which
also makes $f$ identifiable. Our risk bounds imply rates of convergence
for the deconvolution density estimator that are optimal in global
sup-norm loss, without any moment or support restrictions whatsoever,
both in the severely ill-posed case (where linear methods suffice), as
well as in the moderately ill-posed case (where we propose a suitable
thresholding method). To be more precise, given the law $\varphi$ of
the error term and a density $f$ belonging to some Besov body $B(s,L)$
with unknown $s>0$, we devise purely data-driven estimators $\hat f_n$
such that, for every $n \in \mathbb N$,
\[
\sup_{f\in B(s,L)}E\sup_{x\in\mathbb R}|\hat f_n(x)-f(x)|\le r_n(s,\varphi,L),
\]
where $r_n(s, \varphi,L)$ is the minimax rate of convergence in
sup-norm loss over the given Besov body and given the error law
$\varphi$. We also obtain a result of this kind for the case where $f$
is ``supersmooth,'' that is, has an exponentially decaying Fourier
transform. To the best of our knowledge,  the minimax lower bounds
derived in this article are also new.

We should note that the main delicate mathematical point in this work
is to link the $L^2$-based procedure of Fourier inversion to a
pointwise, or even uniform, control of the random fluctuations of the
centered linear density estimator; this problem is already implicit in
the conditions on $F[\varphi]$ and $f$ imposed by Stefanski and Caroll
\cite{SC90}, Fan \cite{F91} and Goldenshluger \cite{Gold1999}, who considered pointwise
loss. Even stronger assumptions were imposed in the nice paper Bissantz
et al. \cite{BDHM07}, wherein  the limiting (extremal-type) distribution of
the uniform deviations over compact sets of certain kernel
deconvolution density estimators for $f$ is derived---this is the only
result that we are aware of in the literature on deconvolution
estimation that deals with sup-norm loss in the moderately ill-posed
case (Stefanski \cite{Stefanski1990} deals only  with the simpler severely ill-posed
case). Our empirical process approach gives results under minimal
conditions and also yields the relevant concentration inequalities that
allow for a satisfactory treatment of adaptation, which the results in
Bissantz et al. \cite{BDHM07} do not address. We should note that applying
empirical process tools in this setting is not at all straightforward:
the usual approach would be to show that certain kernels are of bounded
variation and thus the associated sets of translates and dilates are of
Vapnik--Chervonenkis type (e.g., Nolan and Pollard \cite{NP87}, Einmahl and
Mason \cite{EM00}, Gin\'{e} and Guillou \cite{GG02}), but this does not seem
viable in the deconvolution problem, due to the fact that the bounded
variation norm does not possess a nice Fourier-analytical
characterization. We can, however, solve this problem by combining
recent results on VC properties of functions of quadratic variation in
Gin\'{e} and Nickl \cite{GN09a} with Littlewood--Paley theory and the fact
that wavelet bases are compatible with both the $L^2$- and
$L^\infty$-structure simultaneously; see Lemma \ref{entropy} for this
key result.

Our results can be used to construct confidence bands in the
deconvolution problem and we discuss this in some detail below, as well
as relations to work in \cite{BDHM07,BH}. We suggest a new approach to
nonparametric confidence bands based on Rademacher symmetrization, in a
similar vein as in recent work of Koltchinskii \cite{K06}. While these
confidence bands may be conservative, they allow for an explicit
finite-sample analysis under minimal assumptions.

Let us finally remark that this article also contains new results for
the standard density estimation problem (where $\varphi$ equals Dirac
measure $\delta_0$ at $0$). In this field, our results contribute in
several respects: first, Vapnik--Chervonenkis properties of wavelet
projection kernels have thus far only been derived for Daubechies
wavelets \cite{GN09a} and Battle--Lemari\'{e} wavelets \cite{GN10b},
and the present article achieves the same for wavelets with compactly
supported Fourier transform (e.g., Meyer wavelets). Furthermore, our
main adaptation result, Theorem \ref{result-threshold}, is completely
free of any moment conditions and thus shows, as may have been
suspected, that the moment conditions imposed in Theorem 8 in
\cite{GN09a} are not necessary. Finally, the confidence bands we
suggest can also be used for regular wavelet density estimators and we
are not aware of any other results on global confidence bands in
density estimation, except for the rather technical ones in \cite{gks}.

\section{Main results}

We start with some preliminary definitions and facts. For any Lebesgue
integrable function $h \in L^{1}(\mathbb{R}),$ the Fourier transform
$F[h]$ of $h$ is defined as $F[h](t) =
\int_{\mathbb{R}}h(x)e^{-itx}\,dx$, $t\in \mathbb{R},$ and we use the
natural extension of $F$ to $L^2(\mathbb R)$. We further denote by
$F^{-1}$ the inverse Fourier transform so that $F^{-1}Ff=f$ for $f \in
L^2(\mathbb R)$. The Fourier transform of the density $g$ from
(\ref{deconvolution-model}) is then given by
\begin{equation} \label{conv}
F[g](t) = F[f](t)F[\varphi](t)
\end{equation}
for every $t\in \mathbb{R}$. Another standard property of the Fourier
transform we shall frequently use is its scaling property: for $h\in
L^{1}(\R)$ and $\alpha\in \mathbb R \setminus\{0\},$ the\vspace*{1pt} function
$h_{\alpha}(x):=h(\alpha x)$ has Fourier transform $F[h_{\alpha}](t) =
\alpha^{-1}F[h](\alpha^{-1}t)$.

Let $\phi$ and $\psi$ be, respectively, a scaling function and the
associated wavelet function of a multiresolution analysis. We refer to
\cite{M92,HKPT} for the basic theory of wavelets that we shall use
freely in this article. The dilated and translated scaling and wavelet
functions at resolution level $j$ and scale position $k/2^{j}$ are
defined\vspace*{1.5pt} as $\phi_{jk}(x) = 2^{j/2}\phi(2^{j}x - k)$, $\psi_{jk}(x) =
2^{j/2}\psi(2^{j}x - k)$, $j,k\in\mathbb{Z}$. Now, denote  by
$\langle\cdot,\cdot \rangle$ the inner product in the Hilbert space
$L^{2}(\R)$. The density $f$ can be formally expanded into its wavelet
series
\[
f=\sum_{k\in\mathbb{Z}}\alpha_{jk}(f)\phi_{jk}+\sum_{l=j}^\infty\sum_{k\in\mathbb Z}\beta_{lk}(f)\psi_{lk},
\]
where the coefficients are given by $\alpha_{jk}(f)=\langle f,\phi_{jk}\rangle$, $\beta_{lk}(f)=\langle
f,\psi_{lk}\rangle$,\break
\mbox{$l,j,k\in\mathbb{Z}$}. As is well known, the regularity properties of
a function $f$ can be measured by the decay of their wavelet
coefficients. We define Besov spaces as follows.

\begin{definition} \label{besov}
Let $1\leq p,q\leq\infty$, $s>0$ or let $s=0$ and $q=1$. Let $\phi$
and $\psi$ be the Meyer scaling function and mother wavelet,
respectively (see, e.g., Section 2 of \cite{PV99} for a definition).
The Besov space $B_{pq}^{s}( \R)$ is defined as the set of functions
\[
\Biggl\{f\in L^{p}(\R)\dvtx\|f\|_{s,p,q}=\|\alpha_{0\cdot}\|_{p}+\Biggl(\sum_{l=0}^{\infty}\bigl(2^{l(s+1/2-1/p)}\bigl\|\beta_{l(\cdot)}(f)\bigr\|_{p}\bigr)^{q}\Biggr)^{1/q} <\infty\Biggr\},
\]
where $\|\cdot\|_p$ are the norms of the sequence spaces
$\ell^p(\mathbb Z)$, and with the usual modification in the case
$q=\infty$. Moreover, for any $L>0$, the Besov ball of radius $L$ is defined as
$B(s,p,q,L)=\{f\in L^{p}(\R)\dvtx\|f\|_{s,p,q}\le L\}$.
\end{definition}

\subsection{Minimax lower bounds over Besov bodies}

Before we construct explicit estimators for the density $f$ of $X$ in
the deconvolution model (\ref{deconvolution-model}), we derive a result
that gives a benchmark for the best performance of \textit{any}
estimator $\tilde f_n$. More precisely, we derive lower bounds for the
minimax rate of convergence of $\tilde f_n-f$ in sup-norm loss,
uniformly over Besov bodies of densities $f$ under various assumptions
on the error law $\varphi$. We will  subsequently show that these lower
bounds can be attained by certain wavelet-based estimators and are thus
optimal.

To this end, define the minimax $L^{\infty}$-risk over the H\"older
class $B(s,L):=B(s,\infty,\infty,L)\cap\{f\dvtx\mathbb R\to[0,\infty),\int_{\R}f(x)\,dx=1\}$ as
\begin{equation}\label{minmax}
R_{n}(B(s,L))=\inf_{\tilde{f}_{n}}\sup_{f\in B(s,L)}E\sup_{x\in\mathbb R}|\tilde{f}_{n}(x)-f(x)|,
\end{equation}
where the infimum is taken over all possible estimators
$\tilde{f}_{n}$. Note that an estimator in the deconvolution problem
means any measurable function of a sample $Y_1,\ldots,Y_{n}$ from
density $f\ast\varphi$ that takes values in the space of bounded
functions on $\mathbb R$.

We shall make the following assumption on $F[\varphi]$ to establish the
lower bounds.

\begin{condition}\label{lb-cond}
There exist constants $C,C'>0$, $w,w'\in \mathbb R$ and $t_1,c_0\geq 0$
such that $F[\varphi](t)$ is differentiable for every $t$ satisfying
$|t|>t_1$ and
\[
|F[\varphi](t)|\leq C(1+t^2)^{-w/2}e^{-c_0|t|^{\alpha}},
\]
as well as
\[
|(F[\varphi])'(t)|\leq C'(1+t^2)^{-w'/2}e^{-c_0|t|^{\alpha}}.
\]
\end{condition}

This condition is weaker than the standard ones employed in
deconvolution problems to establish lower bounds
(cf.~\cite{F91,BT08II}), where an additional condition is imposed on
the second derivative of $F[\varphi]$. It covers the usual candidates
for $\varphi,$ including the case $\varphi=\delta_0$ which
corresponds to classical density estimation ($w=c_0=0$).

The following theorem distinguishes the ``moderately ill-posed'' case,
where $F[\varphi]$ decays only polynomially, from the ``severely
ill-posed'' case, where $F[\varphi]$ decays exponentially fast, and
shows that the optimal rates of estimation in the sup-norm depend both
on the smoothness of $f$ and the decay of $F[\varphi]$.

\begin{theorem}\label{theo-minmax}
Let Condition \textup{\ref{lb-cond}} be satisfied. Then, for any $s,L>0,$ there
exists a constant $c:=c(s,L, C,C', \alpha, w, w', c_0)>0$ such that for
every $n\ge 2,$ we have
\[
 R_n(B(s,L))\ge c
 \cases{
\displaystyle\biggl(\frac{1}{\log n}\biggr)^{s/\alpha},&\quad if $c_0 >0,$\cr
\displaystyle\biggl(\frac{\log n}{n}\biggr)^{s/(2s+2w+1)},&\quad if $c_0 =0$ and $w'\geq w\geq 0$.}
\]
\end{theorem}

One may be interested in replacing the H\"{o}lder class $B(s,L)$ by a
more general Besov body, $B(r,p,q,L)$, $r>1/p,$ of densities. It
follows from the proof of Theorem \ref{theo-minmax} that the minimax
rate over $B(r,p,q,L)$ equals the one for $B(s,L)$ with $s=r-1/p$ and
the Sobolev embedding $B^r_{pq}(\mathbb R)\subset B^s_{\infty\infty}(\mathbb R)$ will imply that our upper risk bounds derived in
the following sections attain this rate. We thus restrict ourselves to
$B(s,L)$ without loss of generality.

\subsection{Uniform fluctuations of wavelet deconvolution estimators} \label{linsup}
\subsubsection{The linear wavelet deconvolution estimator}

Recall the model (\ref{deconvolution-model}). We now show, following
\cite{PV99}, how one can estimate $f$ from a sample of $P$ by
``deconvolving'' $P$ or, rather, a suitable approximation of it, on a
wavelet basis $\phi, \psi$ that satisfies the following condition.

\begin{condition} \label{wav}
Assume $\phi,\psi\in L^p(\mathbb R)$ for every $1\le p\le\infty$,
and for some $0<a'<a,$ we have $\operatorname{supp}(F[\phi])\subset[-a,a]$,
as well as $\operatorname{supp}(F[\psi])\subset[-a,-a]\setminus[-a',a']$.
Assume, further, that
\begin{equation}\label{abs}
\qquad c(\phi):=\sup_{x\in\mathbb R}\sum_k|\phi(x-k)|<\infty,\qquad
c(\psi):=\sup_{x\in\mathbb R}\sum_k|\psi(x-k)|<\infty.
\end{equation}
\end{condition}

This condition is satisfied for Meyer wavelets with $a=8\pi/3$ and
$a'=2\pi/3$ (these choices are not optimal, but feasible)---see, for
instance, Section 2 in \cite{PV99}---but other band-limited wavelet
bases are also admissible.

If $K(y,x):=\sum_{k\in \mathbb{Z}}\phi(y-k)\phi(x-k),$ then\vspace*{1.5pt} the
functions $K_j(y,x):=2^jK(2^jy,\break2^jx)$, $j\in\mathbb{N}$, are the
kernels of the orthogonal projections of $L^{2}(\R)$ onto the closed
subspaces $V_j\subset L^2(\mathbb R)$ spanned by $\{\phi_{jk}\dvtx k\in\mathbb Z\}$. We write, for $x\in\mathbb R$, $j\ge 0$ possibly
real-valued,
\begin{eqnarray*}
K_{j}(f)(x)&=&\sum_{k\in\mathbb{Z}}2^{j}\phi(2^{j}x-k)\int_{\R}\phi(2^{j}y-k)f(y)\,dy=\int_{\R}K_{j}(x,y)f(y)\,dy,
\end{eqnarray*}
where the second equality holds pointwise, in view of (\ref{abs}).

Suppose the Fourier transform of the error law $\varphi$ satisfies
$|F[\varphi]| > 0$ on $\operatorname{supp}(F[\phi](2^{-j}(\cdot)))$. We then
have, from Plancherel's theorem, that
\begin{eqnarray} \label{dec}
K_j(f)(x)
&=&
2^j\sum_k\phi(2^jx-k)\int_\mathbb R\phi(2^jy-k)f(y)\,dy\nonumber
\\
&=&
\sum_k\phi(2^jx-k)\frac{1}{2\pi}\int_\mathbb R\overline{F[\phi_{0k}](2^{-j}t)}F[f](t)\,dt\nonumber
\\[-8pt]\\[-8pt]
&=&
\sum_k\phi(2^jx-k)\frac{1}{2\pi}\int_\mathbb R\overline{F[\phi_{0k}](2^{-j}t)}F[g](t)(F[\varphi](t))^{-1}\,dt\nonumber
\\
&=&
2^j\sum_k\phi(2^jx-k)\int_\mathbb R\tilde\phi_{jk}(y)g(y)\,dy=\int_\mathbb R K_j^*(x,y)g(y)\,dy,\nonumber
\end{eqnarray}
where the (nonsymmetric) kernel $K_j^*$ is given by
\[
K_j^*(x,y)=2^j\sum_{k\in\mathbb Z}\phi(2^jx-k)\tilde\phi_{jk}(y)
\]
with
\begin{equation}\label{tilde}
\qquad\tilde\phi_{jk}(x)
=
F^{-1}\biggl[\frac{F[\phi_{0k}](2^{-j}\cdot)}{\overline{2^jF[\varphi]}}\biggr](x)
=
\phi_{0k}(2^j\cdot)\ast F^{-1}\biggl[\frac{1_{[-2^ja,2^ja]}}{\overline{F[\varphi]}}\biggr](x).
\end{equation}
We should note that Young's inequality for convolutions implies, for
fixed $j$, that $\|\tilde\phi_{jk}\|_\infty<\infty$, and then also
$\|K_j^*\|_\infty<\infty$, which justifies the above operations.

Since we have a sample $Y_1,\ldots,Y_n$ from the density $g$, the identity
(\ref{dec}) suggests a natural estimator of $f$, namely the wavelet
deconvolution density estimator
\begin{equation}\label{linear-estimator}
f_n(x,j)=\frac{1}{n}\sum_{m=1}^n K^*_j(x,Y_m),\qquad x\in\mathbb R,\ j\ge 0.
\end{equation}

\subsubsection{Uniform moment and exponential bounds for the fluctuations of $f_n-Ef_n$}

We start with some results for the uniform deviations
\begin{equation} \label{sup}
\qquad\sup_{x\in\mathbb R}|f_n(x,j)-Ef_n(x,j)|\le c(\phi) 2^j\sup_{k\in\mathbb Z}\Biggl|\frac{1}{n}\sum_{m=1}^n\bigl(\tilde\phi_{jk}(Y_m)-E\tilde\phi_{jk}(Y)\bigr)\Biggr|,
\end{equation}
where the inequality follows from (\ref{abs}). This suggests to\vspace*{1.5pt} study
the empirical process indexed by the class of functions $\mathcal F=\{\tilde\phi_{jk}\dvtx k\in\mathbb Z\}$. In fact, some further scaling
depending on the error distribution $\varphi$ will be useful to obtain
a class with constant envelope.

The rather intricate Fourier-analytical definition of $\tilde\phi_{jk}$ in (\ref{tilde}) makes it difficult to apply standard
results from empirical process theory. What is needed is that $\mathcal F$ be a Vapnik--Chervonenkis (VC-type) class of functions. In the
classical density estimation case (where $F[\varphi] =1$), this follows
from results in Nolan and Pollard \cite{NP87} for translates of a fixed
function of bounded variation. We could, however, not control the
bounded variation norm of $\tilde\phi_{jk}$ for general $\varphi$ in a
way that would be useful, mainly because the bounded variation norm
does not interact well with Fourier transforms. Recent results by
Gin\'{e} and Nickl \cite{GN09a} show that the bounded variation condition in
Nolan and Pollard \cite{NP87} can be replaced by $p$-variation for general
$1 \le p < \infty$, and the case $p=2$, which corresponds to
``quadratic variation,'' can be linked in a more efficient way to
Fourier analysis by using Littlewood--Paley theory.

The following key lemma shows that $\mathcal F$, suitably normalized,
is indeed a VC-type class of functions, under minimal conditions on
$F[\varphi]$. Denote by $N(\varepsilon,\mathcal F, L^2(Q))$ the
$\varepsilon$-covering numbers of a class of functions $\mathcal F$
with respect to~the $L^2(Q)$-distance.

\begin{lemma} \label{entropy}
Suppose that $\phi,\psi$ satisfy Condition \textup{\ref{wav}} and that
$|F[\varphi](t)|>0$ on $[-2^ja,2^ja]$. Define
\begin{equation}\label{del}
\delta_j:=\min_{t\in[-2^ja, 2^ja]}|F[\varphi](t)|
\end{equation}
(which exists and is positive for every $j$ since $\varphi$ is a
probability measure). Then the class $\mathcal H_j = \{\delta_j
\tilde\phi_{jk}\dvtx k\in\mathbb Z\}$, $j\ge 0$, is uniformly
bounded by the constant $U$ and satisfies, for every $0<\varepsilon<A$, $\sup_{Q} N(\varepsilon,\mathcal H_j,L^2(Q))\le(A/\varepsilon )^v$ for finite positive constants $A,v,U$
depending only on $\phi,\psi$, and where the supremum extends over all
probability measures $Q$ on $\mathbb R$.
\end{lemma}

Combining this lemma with moment bounds for empirical processes indexed
by VC-type classes of functions in \cite{EM00,GG02}, as well as with
Talagrand's \cite{talagrand1996} inequality, we obtain the following
result.

\begin{proposition}\label{suprat}
Suppose that $\phi,\psi$ satisfy Condition \textup{\ref{wav}}, that
$|F[\varphi](t)|>0$ on $[-2^ja,2^ja]$, let $\delta_j$ be as in
\textup{(\ref{del})} and define $j'=\max(1,j)$. Let $f_n(x,j)$ be the
deconvolution wavelet density estimator from
\textup{(\ref{linear-estimator})}
and assume that $X$ has a bounded density $f\dvtx\mathbb R\to[0,\infty)$. Then there  exists a constant $L'$, depending only on $\phi,\psi,p$, such that for every $n\ge 1$, every $j\ge 0$ and $1\le p<\infty$,
\[
\Bigl(E\Bigl(\sup_{x\in\mathbb R}|f_n(x,j)-Ef_n(x,j)|\Bigr)^p\Bigr)^{1/p}\le L'\frac{1}{\delta_j}\biggl(G\sqrt{\frac{2^jj'}{n}}+\frac{2^j j'}{n}\biggr),
\]
where $G=\max(\|g\|^{1/2}_\infty,1)$.
In addition, there exists a constant $C,$ depending only on $\phi,\psi$, such that for every $j\geq 0$ and $u>0,$
\begin{eqnarray}\label{lin-tal}
&&\Pr\biggl\{\sup_{x\in\mathbb
R}|f_n(x,j)-Ef_n(x,j)|\ge\frac{C}{\delta_{j}}\biggl(G\sqrt{(1+u)\frac{2^jj'}{n}}+(1+u)\frac{2^jj'}{n}\biggr)\biggr\}\nonumber
\\[-8pt]\\[-8pt]
&&\qquad\le
e^{-(1+u)j'}.\nonumber
\end{eqnarray}
\end{proposition}

The constant $C$ is unspecified here, although it could be computed
explicitly. Obtaining realistic constants is an intricate matter, but
one can use symmetrization techniques to circumvent this problem; see
Proposition \ref{rad} below.

\subsubsection{Uniform fluctuations of the empirical wavelet coefficients}

The techniques from the previous section allow us to establish similar
uniform estimates for the deviations of the empirical wavelet
deconvolution coefficients $\hat{\beta}_{lk}$ from their means. Such
results are particularly interesting for nonlinear thresholding
procedures that we shall study below.

We have, for $\psi$ satisfying Condition \ref{wav},
\begin{eqnarray*}
\beta_{lk}(f)
&=&
2^{l/2}\int_{\R}\psi(2^lx-k) f(x)\,dx
=
\frac{2^{l/2}}{2\pi}\int_{\R}2^{-l}\frac{\overline{F[\psi_{0k}](2^{-l}\cdot)}}{F[\varphi]}(t)F[g](t)\,dt
\\
&=&
2^{l/2}\int_{\R}F^{-1}\biggl[2^{-l}\frac{F[\psi_{0k}](2^{-l}\cdot)}{\overline{F[\varphi]}}\biggr](x)g(x)\,dx
=:
2^{l/2}\int_{\R}\tilde\psi_{lk}(x)g(x)\,dx.
\end{eqnarray*}
A natural unbiased estimator of $\beta_{lk}\equiv \beta_{lk}(f)$ is
therefore
\begin{equation}
\hat{\beta}_{lk}(f)=\frac{2^{l/2}}{n}\sum_{m=1}^{n}\tilde{\psi}_{lk}(Y_{m})
\end{equation}
and the object of interest in this subsection is the random variable
$\sup_{k\in\mathbb Z}|\hat\beta_{lk}-\beta_{lk}|$. We should note
that for wavelets satisfying Condition \ref{wav} (e.g., Meyer
wavelets), and even if $g$ has compact support, the last supremum is
over an infinite set, so empirical process techniques are particularly
useful. Lemma \ref{entropy} and Proposition \ref{suprat} have the
following analogs for $\tilde\psi$.

\begin{lemma} \label{entropy1}
Suppose that $\phi,\psi$ satisfy Condition \textup{\ref{wav}}, that
$|F[\varphi](t)|>0$ on $[-2^la,2^la]$ and let $\delta_l$ be as in
\textup{(\ref{del})}. Then the class $\mathcal D_l=\{\delta_l\tilde\psi_{lk}\dvtx k\in\mathbb Z\}$, $l\ge 0,$ is uniformly bounded by a fixed
constant $U$ and satisfies, for every $0<\varepsilon<A$, $\sup_{Q}N(\varepsilon,\mathcal D_l,L^2(Q))\le(A/\varepsilon)^v$ for constants $U,A,v$ depending only on $\phi,\psi$.
\end{lemma}

\begin{proposition}\label{suprat-beta}
Suppose that $\phi, \psi$ satisfy Condition \textup{\ref{wav}}, that
$|F[\varphi](t)|>0$ on $[-2^la,2^la]$, let $\delta_l$ be as in
\textup{(\ref{del})} and define $l'=\max(l,1)$. Assume that $X$ has a bounded
density $f\dvtx\mathbb R\to [0, \infty)$. Then, for every $n\ge 1$, for
every $l\geq 0$ and $1 \le p < \infty$, we have
\[
\Bigl(E\sup_{k\in\mathbb{Z}}|\hat{\beta}_{lk}-\beta_{lk}|^{p}\Bigr)^{1/p}
\leq
L''\frac{1}{\delta_{l}}\biggl(G\sqrt{\frac{l'}{n}}+\frac{2^{l/2}l'}{n}\biggr),
\]
where $L''>0$ depends only on $p,\phi,\psi$ and where $G$ is as in
Proposition \ref{suprat}. In addition, there exists a constant $D,$
depending only on $\phi,\psi,$ such that for every $l\geq 0$ and
$u>0$,
\begin{equation}\label{beta-tal}
\qquad\Pr\biggl\{\sup_{k\in\mathbb{Z}}|\hat{\beta}_{lk}-\beta_{lk}|
\ge
\frac{D}{\delta_{l}}\biggl(G\sqrt{(1+u)\frac{l'}{n}}+(1+u)\frac{2^{l/2}l'}{n}\biggr)\biggr\}
\le
e^{-(1+u)l'}.\hspace*{-12pt}
\end{equation}
\end{proposition}

\subsection{Optimal estimation over H\"older classes}

We now show how the risk bounds from the previous section imply optimal
rates of convergence for densities $f\in B^s_{\infty\infty}(\mathbb
R)$ in the deconvolution problem, under the standard decay conditions
on $F[\varphi]$ from the inverse problem literature.

We first consider the case where the error law $\varphi$ decays
exponentially fast. In this ``severely ill-posed'' case, one can find a
universal choice of $j$ for which the linear estimator attains the
exact minimax rate, even without having to know the value $s$.

\begin{theorem}\label{uni-band}
Suppose that $\phi, \psi$ satisfy Condition \textup{\ref{wav}} and assume that\break $|F[\varphi](t)| \ge Ce^{-c_0|t|^\alpha}$ for every $t \in \mathbb R $
and some $C,c_0,\alpha>0$. Let $f_{n}(\cdot, j_{n})$ be the estimator
defined in \textup{(\ref{linear-estimator})}, where
$j_{n}=\frac{1}{\alpha}\log_{2}(\nu\log n)$ for some $\nu$ satisfying
$c_0 a^\alpha\nu <1/2$. Then there exists a constant $L'''$, depending
only on $s,L,\phi,\psi,c_{0},C,\alpha,\nu$, such that for every
$n\ge 2$, we have
\[
\sup_{f \in B(s,L)}E\sup_{x\in\mathbb R}|f_{n}(x,j_{n})-f(x)|
\leq
L'''\biggl(\frac{1}{\log n}\biggr)^{s/\alpha}.
\]
\end{theorem}

We now turn to the case where $F[\varphi]$ decays polynomially, the
so-called ``moderately ill-posed'' case. Here, the linear estimator
$f_n$ is only minimax optimal if one knows the value of $s$.\vspace*{-1pt}

\begin{theorem}\label{mod-ill-theo}
Suppose that $\phi,\psi$ satisfy Condition \textup{\ref{wav}}\vspace*{1.5pt} and assume that $|F[\varphi](t)|\ge C(1+|t|^{2})^{-w/2}$ for every $t\in\mathbb R$ and some $C>0$, $w\ge 0$. Let\vspace*{1.5pt} $f_{n}(\cdot,j_{n})$ be the
estimator defined in \textup{(\ref{linear-estimator})} with $j=j_{n}$ satisfying
$2^{j_{n}}\simeq\break (n/\log n)^{1/(2s+2w+1)}$. Then there
exists a constant $C'$, depending only on $s,L,\phi,\psi,C,w$, such
that for every $n\ge 2$, we have
\[
\sup_{f\in B(s,L)}E\sup_{x\in\mathbb R}|f_{n}(x,j_{n})-f(x)|\leq C'\biggl(\frac{\log n}{n}\biggr)^{s/(2s+2w+1)}.
\]
\end{theorem}

The question arises as to whether we can achieve this rate of
convergence without having to know the value of $s$ in our choice of
$j_n$ so that we can adapt to the unknown smoothness $s$ of $f$. This
can be done using the wavelet thresholding deconvolution estimator
proposed in Johnstone et al. \cite{JKPR04} in the
periodic setting, defined as follows: for $j_{1}$ positive integers, to
be specified below, the hard thresholding estimator equals
\begin{equation}\label{est-threshold}
f_n^{T}(x)
=
f_n(x,0)+\sum_{l=0}^{j_{1}-1}\sum_{k}\hat{\beta}_{lk}1_{|\hat{\beta}_{lk}|>\tau}\psi_{lk}(x),
\end{equation}
where $\hat\beta_{lk}$ was introduced in Section \ref{linsup}. The
threshold $\tau$ is chosen such that $\tau =\tau(n,l,w,\kappa)=\kappa 2^{wl}\sqrt{(\log n)/n},$ where $\kappa=G\kappa'$, with $G$
from Proposition \ref{suprat} and $\kappa'$ a ``large enough'' constant
that depends only on $w,C,\phi,\psi$. If $G$ is unknown, then it can
be replaced by an estimate, as in \cite{GN10b}.

\begin{theorem}\label{result-threshold}
Suppose that $\phi,\psi$ satisfy Condition \textup{\ref{wav}}. Suppose that
$\varphi$ is such that $|F[\varphi](t)|\geq C(1+|t|^{2})^{-w/2}$ for every $t\in\mathbb R$ and some
$C>0$, $w\ge 0$. Let $f_{n}^T$ be the thresholded estimator in
\textup{(\ref{est-threshold})} with
\begin{eqnarray*}
\biggl(\frac{n}{\log n}\biggr)^{1/(2w+1)}
&\le&
2^{j_{1}}\leq 2\biggl(\frac{n}{\log n}\biggr)^{1/(2w+1)},\qquad j_{1}>0.
\end{eqnarray*}
We then have, for every $n \ge 2$ and every $s>0$, that
\begin{equation} \label{maint}
\sup_{f\in B(s,L)}E\sup_{x\in\mathbb R}|f_{n}^{T}(x) -f(x)|
\leq
D\biggl(\frac{\log n}{n}\biggr)^{s/(2w+2s+1)},
\end{equation}
where $D>0$ depends only on $s,L,\phi,\psi,w,C$.
\end{theorem}

\subsection{Extensions and applications} \label{app}

\subsubsection{Estimation of a supersmooth density}

In the last sections, we established the minimax rate of estimation of
a density in $B^{s}_{\infty\infty}(\mathbb R)$ for the sup-norm error,
both in the moderately and severely ill-posed cases, and constructed
estimators that attain this rate. It was pointed out in \cite{PV99}
for the $L^{2}$-error that the linear and thresholded estimators attain
faster rates of convergence if we consider classes of supersmooth
densities instead of the usual Besov spaces. In this section, we
investigate this phenomenon for the sup-norm error. We show that the
minimax rate of convergence for the sup-norm is the same as that
obtained for the $L^{2}$-error up to an additional $\sqrt{\log\log n}$
factor and that wavelet estimators can attain this rate. For
simplicity, and to highlight the main ideas, we only consider the
nonadaptive case.

Assume that $f$ belongs to the class of supersmooth densities,
\[
\mathcal A_{\tilde{c}_0,s}(L)
=
\biggl\{f\dvtx\mathbb{R}\rightarrow [0,\infty),\int_{\mathbb{R}}f
=
1,\int_{\mathbb R}|F[f](t)|^{2}\exp(2\tilde{c}_0|t|^{s})\,dt\le 2\pi L\biggr\},
\]
where $\tilde{c}_0,s,L> 0$. In the moderately ill-posed case, we have
the following result.

\begin{corollary} \label{supsm}
Let $\phi,\psi$ satisfy Condition \textup{\ref{wav}}. Assume that $f\in \mathcal
A_{\tilde{c}_0,s}(L)$ for some $\tilde{c}_0,s,L>0$ and that
$|F[\varphi](t)|\ge C(1+|t|^2)^{-w/2}$ for every $t \in \mathbb
R$ and some $C>0$, $w\ge 0$. Let $f_{n}(\cdot,j_{n})$ be the estimator
defined in \textup{(\ref{linear-estimator})} with $j=j_{n}$ satisfying
\[
2^{j_{n}}
=
\biggl(\frac{1}{2(a')^{s}\tilde{c}_0}\log n\biggr)^{1/s}.
\]
Then there exists a constant $C',$ depending
only on $\phi,\psi, \tilde{c}_0,s, L, C,w,$ such that for every $n\ge
3,$ we have
\[
\sup_{f\in\mathcal A_{\tilde c_0,s}(L)}E\sup_{x\in\mathbb R}|f_{n}(x,j_{n})-f(x)|
\leq
C'\biggl(\frac{\log\log n}{n}\biggr)^{1/2}(\log n)^{(w+1/2)/s}.
\]
\end{corollary}

The rates we obtained for the sup-norm error are similar to those
obtained by \cite{BT08I,BT08II} and \cite{PV99} for the $L^{2}$-error,
up to the presence of the additional factor $\sqrt{\log\log n}$. This
additional factor can be heuristically\vspace*{1pt} explained by the presence of the
quantity $\sqrt{j}$ in the deviation term
$\delta_j^{-1}(2^{j}j/n)^{1/2}$ derived in\vspace*{1pt} Proposition \ref{suprat}.
The next theorem implies that this $\sqrt{\log\log n}$ factor is
indeed necessary.

\begin{theorem} \label{minm2}
Fix $0<s\le 1$ and $\tilde{c}_0,L>0$. Assume that $\varphi$ satisfies
Condition \textup{\ref{lb-cond}} with $c_0=0$ and $w'\geq w\geq 0$. There
then exists a positive constant $c:=c(s,\tilde{c}_0,L,C,C',w,w')$
such that
\[
\inf_{\tilde{f}_{n}}\sup_{f\in\mathcal{A}_{\tilde{c}_0,s}(L)}E\sup_{x\in\mathbb R}|\tilde{f}_{n}(x)-f(x)|
\ge
c\biggl(\frac{\log\log n}{n}\biggr)^{1/2}(\log n)^{(w+1/2)/s}.
\]
\end{theorem}

We can also obtain a faster rate of convergence in the severely
ill-posed case for supersmooth densities, balancing the bias bound from
Proposition \ref{smooth-bias} below with the variance bound from
Proposition \ref{suprat} above. We can then obtain similar results as
in \cite{BT08I,BT08II}, with additional logarithmic terms in the rate
of convergence, due to the fact that we consider sup-norm loss instead
of $L^2$-loss.

\subsubsection{Confidence bands}

One of the main statistical challenges in the nonparametric
deconvolution problem is the construction of confidence bands for $f$
(cf.~\cite{BDHM07,BH}). In \cite{BDHM07}, the exact uniform (over
compact subsets of $\mathbb R$) limit distribution of certain linear
kernel-based deconvolution estimators for $f$ is derived, assuming that
$f$ satisfies $\int_\mathbb R|F[f](u)||u|^r\,du<\infty$ for $r>0$ and
that $g$ is once differentiable with bounded derivative, and if the
Fourier transform of the error variable decays exactly like a
polynomial, that is, $|F[\varphi](t)|\asymp C|t|^{-w}$ for some
$C>0$, $w\ge 0$. If the underlying smoothness $r$ of $f$ is known,
then these results can be used to construct asymptotic confidence bands
for $f$ that shrink at certain rates of convergence.

We suggest here an alternative approach to confidence bands in the
nonparametric deconvolution problem. Instead of extreme value theory,
we use concentration inequalities and Rademacher processes. This allows
for almost assumption-free results and has the advantage that the
confidence band can be shown to be valid on the whole real line and for
every sample size $n$. On the downside, these bands are likely to be
too conservative in the limit.

One fundamental problem of using concentration inequalities (as in
Proposition \ref{suprat}) in practice is that often, no reasonable
values for the leading constant $C$ are available. To circumvent this
problem, we use here an idea that goes back to Koltchinskii \cite{K01,K06} and Bartlett, Boucheron and Lugosi \cite{BBL01}; see also Gin\'{e} and
Nickl \cite{GN10b}, where this approach was introduced in density
estimation. Define a Rademacher process
 and the associated supremum,
\[
\Biggl\{\frac{1}{n}\sum_{m=1}^{n}\varepsilon_m K^*_j(x,Y_m)\Biggr\}_{x\in\mathbb R},
\qquad R_n(j):=\sup_{x\in\mathbb R}\Biggl|\frac{1}{n}\sum_{m=1}^{n}\varepsilon_m K^*_j(x, Y_m)\Biggr|,
\]
with $(\varepsilon_m)_{m=1}^n$ an i.i.d.~Rademacher sequence, independent of
the $Y_m$'s (and defined on a large product probability space). $R_n$
can be computed in practice by first simulating $n$ i.i.d.~random
signs, applying these signs to the summands $K^*_j(x,Y_m)$ of the
wavelet deconvolution density estimator (\ref{linear-estimator}) and
maximizing the resulting function. Similarly, one can consider
$E_\varepsilon R_n(j)$, the expectation of $R_n(j)$ with respect to the
Rademacher variables only, which is a stochastically more stable
quantity.

We shall use the fact that this is the supremum of a centered process
which can be shown to concentrate around
$2E\|f_n(\cdot,j)-Ef_n(\cdot,j)\|_\infty$. To describe the
concentration property, recall $\delta_j$ from (\ref{del}) and define
the random variable
\begin{equation}\label{sigr}
\quad\qquad\sigma^{R}(n,j,z)
=
6R_n(j)+\frac{D_1}{\delta_j}\sqrt{\frac{2^{j}\|g\|_\infty(z+\log 2)}{n}}
+
\frac{D_2}{\delta_j}\frac{2^{j}(z+\log2)}{n},
\end{equation}
where $D_1=10 c(\phi)\|\phi\|_1\sqrt{a/\pi}\le 5.7 c(\phi)\|\phi\|_1\sqrt a $,
$D_2 = 44 c(\phi)\sqrt{a/2\pi^2}\le\break 11 c(\phi)\sqrt a$
and $c(\phi)$ as in (\ref{abs}). If $\|g\|_\infty$ is unknown,
it can be replaced by $\|f_{n}(\cdot,j_n)\|_\infty$ in practice so that
$\sigma^R$ is completely\vspace*{1pt} data-driven. We start with a confidence band
$\bar C_n=[f_n(\cdot,j)-\sigma^R(n,j,z),f_n(\cdot,j)+\sigma^R(n,j,z)]$ for the mean $Ef_n$ of $f_n$.

\begin{proposition}\label{rad}
Let $f_n(x,j)$ be the estimator from \textup{(\ref{linear-estimator})} and
suppose that $|F[\varphi]|>0$ on $[-2^ja,2^ja]$. Assume that $X$ has a
bounded density $f\dvtx \mathbb R\to[0,\infty)$. We then have, for every
$n\ge 1$, every $j\in\mathbb N$ and every $z>0,$ that
\[
\Pr\Bigl\{\sup_{x\in\mathbb R}|f_n(x,j)-Ef_n(x,j)|\ge\sigma^R(n,j,z)\Bigr\}\le e^{-z}.
\]
Moreover, the band $\bar C_n$ has expected diameter
\[
2E\sigma^R(n,j,z)\le C\delta^{-1}_j\biggl(\sqrt{\frac{2^jj}{n}}+\frac{2^jj}{n}\biggr)
\]
for every $z>0$, every $n\in\mathbb N$, every
$j\ge 1$ and some constant $C$ depending only on $\|g\|_\infty,\phi,\psi,z$.
\end{proposition}

Proposition \ref{rad} still holds true when $R_n(j)$ is replaced by
$E_{\varepsilon}R_n(j)$, the expectation of $R_n(j)$ with respect
to~the Rademacher variables only. (This follows from combining the
proof of Proposition \ref{rad} with the arguments in the proof of
Proposition 2 in \cite{GN10b}.)

We did not try to optimize the constants in the choice of $\sigma^R$
and they are likely to be suboptimal, as they depend on the constants
in the lower-deviation version of Talagrand's inequality, where sharp
constants are not yet known. A ``practical'' choice may be to replace
the $6$ in front of $R_n$ by $4$ and to ignore the third ``Poissonian''
term in (\ref{sigr}).

We again emphasize  that we simply need $|F[\varphi](t)|$ to be bounded
from below on the fixed interval $[-2^ja,2^ja]$ for our results to hold
and we do not need any support or moment assumptions on $f$. In
particular, this nonasymptotic result can even be used in principle
when $F[\varphi]$ equals zero eventually, by choosing $j$ small enough.

If $f\in B^s_{\infty\infty}(\mathbb R)$, with $s$ known, then the
last proposition can be readily applied for the construction of
confidence bands $C_n$ for the unknown density $f$ using undersmoothing
(just as in \cite{BDHM07}) and these bands can be shown to shrink at
the optimal rate of convergence depending on the smoothness of $f$. We
do not detail this here, nor do we address the more difficult problem
of \textit{adaptive} confidence bands: using Proposition \ref{rad},
such results can be obtained in the same way as in the case of
density estimation considered in \cite{GN10a}.

Instead, and for sake of illustration, let us construct a nonasymptotic
confidence band in the supersmooth case $f\in\mathcal A_{\tilde{c}_0,s}(L)$, $s,\tilde c_0$ known, with moderately ill-posed
error distribution.

\begin{corollary}\label{band}
Let $f,\varphi,f_{n}(\cdot,j_{n})$ and $j_n$ be as in Corollary \textup{\ref{supsm}}. Let $\sigma^R(n,j,z)$ be as in \textup{(\ref{sigr})} above and define the confidence band
\[
C_n(x,z)=[f_n(x,j_n)\pm(1+\delta)\sigma^R(n,j_n,z)],\qquad x\in\mathbb R,
\]
where $\delta$ is any positive real number. Then, for every $z>0$ and every $n\in\mathbb N$,
\[
\Pr\{f(x)\in C_n(x,z)\ \forall x\in\mathbb R\}\ge 1-e^{-z}-v_n,
\]
where [$c'''\equiv c'''(\phi,\psi,\tilde c_0,s)$, as in Proposition \textup{\ref{smooth-bias}}]
\[
v_n\equiv\Pr\biggl\{\sigma^R(n,j_n,z)\le\frac{c'''}{\delta}\sqrt{L\frac{(\log n)^{(1-s)/s}}{n}}\biggr\}
\]
satisfies $v_n\to 0$ as $n\to\infty$.

Moreover, if $|C_n(z)|$ is the maximal diameter of $C_n(x,z),$ then
\[
E|C_n(z)|\le C\biggl(\frac{\log\log n}{n}\biggr)^{1/2}(\log n)^{(w+1/2)/s},
\]
where $C$ depends on $\tilde{c}_0,s,L,\delta,z,\|g\|_\infty$.
\end{corollary}

Since $\lim v_n =0,$ this confidence band has asymptotic coverage for
$\delta>0$ arbitrary, but more is true: $v_n$ equals zero from some $n$
onward and one can, in principle, even obtain coverage for every fixed
sample size $n$ by choosing $\delta$ in dependence of $L$ (and of the
constants that define $\sigma^R$).

\section{Proofs}

\subsection{\texorpdfstring{Proof of Theorem
\protect\ref{theo-minmax}}{Proof of Theorem 1}}

Our proof adapts to the present situation standard lower bound
techniques as in \cite{BT08II,F91,PS09}. We recall that the
Kullback--Leibler divergence between two distributions $P$ and $Q$ is
defined by
\[
K(P|Q)
=
\cases{\displaystyle\int\log\biggl(\frac{dP}{dQ}\biggr)\,dP,&\quad if $P\ll Q,$\cr
+\infty,&\quad elsewhere.}
\]
To establish lower bounds for the minimax risk (\ref{minmax}), we use
the following lemma (see Theorem 2.5 on page 99 of
\cite{Tsy2009})---actually, an adaptation of it---to the deconvolution
problem at hand.

\begin{lemma}\label{lemlb}
Let $d$ be a metric on $B(s,L)$. Let $r_n$ be a sequence of positive
real numbers and let $\mathcal C\subset B(s,L)$ be a finite set of
probability densities such that $\operatorname{card}(\mathcal C)\ge 2$ and $\forall f,g\in\mathcal C$, $f\neq g\Rightarrow d(f,g)\ge 4 r_n>0$.
Further, let $\varphi$ be a fixed probability measure and let $P_{f\ast\varphi}^n$ be the product probability measure corresponding to a
sample of size $n$ from the law $f\ast\varphi$, $f\in\mathcal C$,
and assume that the KL divergences satisfy, for every $f\in\mathcal C$ and some $f_0\in\mathcal C$,
\[
K(P_{f\ast\varphi}^n |P_{f_0\ast\varphi}^n)\le\frac{1}{16}\log(\operatorname{card}(\mathcal C)).
\]
Then,
\[
\inf_{\hat{f}_n}\sup_{f\in\mathcal C} E d(\hat f_n,f)\ge c_1 r_n,
\]
where $\inf_{\hat{f}_n}$ denotes the infimum over all estimators based\vspace*{1pt}
on a sample of size $n$ from the density $f\ast\varphi$ and where
$c_1>0$ is a constant that depends only on $s,L$.
\end{lemma}

We use this lemma to prove Theorem \ref{theo-minmax}. Let $\psi$ be the
Meyer wavelet. Fix $s,L>0$ and let $j\in\mathbb N$ be arbitrary (to
be chosen later). Define the set of functions $\mathcal C=\{f_k,k=0,\ldots,2^{j}-1\}$ as follows: consider the standard Cauchy density
$p(x)=1/\pi(1+x^2)$, set $f_0(x)=\frac{1}{\eta}p(\frac{x}{\eta})$
for $\eta>0$ and for any $k=1,\ldots,2^j-1$, $f_k(x)=f_{0}(x)+c'2^{-j(s+1/2)}\psi_{jk_M},$ where $k_M=Mk$ for some integer $M\ge 1$
specified below. We show that the constants $\eta,c'>0$ can be chosen
such that $f_k$ is a density on $\mathbb R$ and, in fact, belongs to
$B(s,L)$ for every $k=1,\ldots,2^{j}-1$ and every integer~$M$. Clearly,
$f_k$ integrates to 1 since $\psi$ is orthogonal on constants. We next
prove $f_k\in B(s,L)$ for all $k$ and suitable $c',\eta$. First, we
have $\|f_0\|_{s,\infty,\infty}\le\frac{L}{2}$ for $\eta\ge 1$ large
enough and depending only on $s,L,\psi,\phi$, in view of $F[f_0](u)=e^{-\eta|u|}$,
Definition \ref{besov}, $|\beta_{lk}(f_0)|=|(1/2\pi)\int_{\R} e^{-\eta|u|}F[\psi_{lk}](u)|\le 2^{-l/2}\|\psi\|_1 e^{-|2^la'|\eta}$
with $a'=2\pi/3$ and a similar estimate
for $\alpha_k(f_0)$. Thus, we have, for $0<c'\le L/2,$
\[
\|f_k\|_{s,\infty,\infty}
\le
\|f_0\|_{s,\infty,\infty}+\bigl\|c'2^{-j(s+1/2)}\psi_{jk_M}\bigr\|_{s,\infty,\infty}
\le
\frac{L}{2}+c'\le L.
\]
Having chosen $\eta$, we can choose $c'\le L/2$ suitably small but
positive and depending on $\eta$ and $\psi$ so that $f_k>0$ on
$\mathbb R$ for any $k$. This  is easily established by using the fact
that the Meyer wavelet decays faster at infinity than any polynomial
[i.e., the estimate $|\psi(x)|\le C_N/(1+|x|^2)^{N/2}$ for
every $N\in\mathbb N$ and every $x \in \mathbb R$], whereas $f_0(x)$
decays at infinity like $x^{-2}$.

To proceed with the proof, we set $\gamma_j =
c'2^{-j(s+1/2)}$. We first prove the separation property in sup-norm
for the $f_k$'s. For\vspace*{1pt} any distinct $f_{k}, f_{k'}$, we have $\|f_k -
f_{k'}\|_\infty = \gamma_{j}2^{j/2}\|\psi(\cdot -Mk) -
\psi(\cdot-Mk')\|_\infty.$ By definition of the Meyer wavelet, we have,
for any $k\neq k'$,
\begin{eqnarray*}
\sup_{x}|\psi(x-Mk)-\psi(x-Mk')|
&=&
\sup_{x}\bigl|\psi(x)-\psi\bigl(x+M(k-k')\bigr)\bigr|
\\
&\ge&
\|\psi\|_\infty-\bigl|\psi\bigl(x_{\max}+M(k-k')\bigr)\bigr|
\end{eqnarray*}
for some $x_{\max}\in\operatorname{arg\,max}_x|\psi(x)|$. By the decay property of
the Meyer wavelets mentioned above, there exists a numerical constant
$M\ge 1$, large enough but finite, such that for any $x$ satisfying
$|x|\ge M$, we have $|\psi(x_{\max}+x)|\le\|\psi\|_\infty/2$. Thus,
we have, for any $k\neq k'$,
\[
\|f_k-f_{k'}\|_\infty\ge\gamma_j 2^{j/2}\frac{\|\psi\|_\infty}{2}
=
2^{-js}\frac{c'\|\psi\|_\infty}{2}.
\]

We now check the second condition of Lemma \ref{lemlb}. Let
$(Y_{1},\ldots,Y_{n})$ be an i.i.d.~sample with distribution $P^n_{k}$
admitting the density $\prod_{i=1}^{n}(f_{k}\ast \varphi)(y_{i})$ with
respect to~the Lebesgue measure on $\mathbb R^{n}$. Fubini's theorem
and the fact that $\psi$ is orthogonal on constants give, for $k\in\mathbb Z$,
that $\int_{\mathbb R}(\psi_{jk}\ast\varphi)(y)\,dy=0$.
Thus, by definition of the Kullback--Leibler divergence and the
inequality $\log(1+x)\le x$ for $x>-1$, we obtain, for any $k=1,\ldots,2^{j}-1$, that
\begin{eqnarray}\label{interm-minmax-1}
K(P^n_{k}|P^n_{0})
&=&
n\int_{\mathbb R}\log\biggl(\frac{f_{k}\ast\varphi}{f_{0}\ast\varphi}(y)\biggr)(f_{k}\ast\varphi)(y)\,dy\nonumber
\\
&=&
n\int_{\mathbb R}\log\biggl( 1+ \gamma_{j}\frac{\psi_{jk_M}\ast\varphi}{f_0\ast\varphi}(y)\biggr)(f_{k}\ast\varphi)(y)\,dy\nonumber
\\[-8pt]\\[-8pt]
&\le&
n\gamma_{j}\int_{\mathbb R}(\psi_{jk_M}\ast\varphi)(y)\biggl(1+\gamma_j\frac{\psi_{jk_M}\ast\varphi}{f_0\ast\varphi}(y)\biggr)\,dy\nonumber
\\
&\le&
n\gamma_{j}^{2}\int_{\mathbb R}\frac{(\psi_{jk_M}\ast\varphi)^{2}}{f_0\ast\varphi}(y)\,dy.\nonumber
\end{eqnarray}
To proceed, we observe that $f_0$ being Cauchy implies that $(f_{0}\ast\varphi)(y)\ge c_{1}/(1+y^{2})$ for some $c_1>0$ and every
$y\in\mathbb{R}$. This is obviously true for $y$ in any compact set
$[-A,A]$, and for $|y|>A$, it follows from
\[
\liminf_{|y|\to\infty}(1+y^2)f_0\ast\varphi(y)
\ge
\frac{1}{\eta\pi}\int_{\R}\liminf_{|y|\to\infty}\frac{1+y^2}{1+[(y-x)/\eta]^2}\,d\varphi(x)
=
\frac{\eta}{\pi},
\]
in view of Fatou's lemma. Consequently, we have
\[
\int_{\mathbb R}\frac{(\psi_{jk_M}\ast\varphi)^{2}}{f_0\ast\varphi}(y)\,dy
\le
\frac{1}{c_{1}}\int_{\mathbb R}(1+y^{2})(\psi_{jk_M}\ast\varphi)^{2}(y)\,dy.
\]
Let us first consider the quantity $\int_{\mathbb R}(\psi_{jk_M}\ast\varphi)^{2}(y)\,dy$. Plancherel's theorem gives
\begin{eqnarray}\label{minmax-interm2}
\qquad\quad\int_{\mathbb R}(\psi_{jk_M}\ast\varphi)^{2}(y)\,dy
&=&
c_{2}\int_{\mathbb R}|F[\psi_{jk_M}](t)|^{2}|F[\varphi](t)|^{2}\,dt\nonumber
\\[-8pt]\\[-8pt]
&\le&
c_3 2^{-j}\|\psi\|_{1}^{2}\int_{\operatorname{supp}(F[\psi_{jk_M}])}(1+t^{2})^{-w}e^{-2c_0|t|^{\alpha}}\,dt\nonumber
\end{eqnarray}
for some constants $c_{2},c_3>0$ depending\vspace*{1pt} only on $C,\pi$.

For the quantity $\int_{\mathbb R}y^{2}(\psi_{jk_M}\ast\varphi)^{2}(y)\,dy$, we obtain similarly, using in addition the
spectral representation of the differential operator, that
\begin{eqnarray}\label{minmax-interm3}
&&
\int_{\mathbb R}(y\psi_{jk_M}\ast\varphi)^{2}(y)\,dy\nonumber
\\
&&\qquad=
c_{2}\int_{\mathbb R}|(F[\psi_{jk_M}](t)F[\varphi](t))'|^{2}\,dt\nonumber
\\
&&\qquad=
c_2\int_{\mathbb R}\bigl|\bigl(F[\psi_{jk_M}]'(t)F[\varphi](t)+F[\psi_{jk_M}](t)F[\varphi]'(t)\bigr)\bigr|^{2}\,dt
\\
&&\qquad=
c_2\int_{\mathbb R}\bigl|\bigl(2^{-3j/2}F[\psi_{0k_M}]'(2^{-j}t)F[\varphi](t)+F[\psi_{jk_M}](t)F[\varphi]'(t)\bigr)\bigr|^{2}\,dt\nonumber
\\
&&\qquad\le
2c_4 2^{-3j}\biggl(\int_{\R}|x\psi(x)|\,dx\biggr)^{2}\int_{\operatorname{supp}(F[\psi_{jk_M}])}(1+t^{2})^{-w}e^{-2c_0|t|^{\alpha}}\,dt\nonumber
\\
&&\qquad\quad{}+
2c_42^{-j}\|\psi\|_{1}^{2}\int_{\operatorname{supp}(F[\psi_{jk_M}])}(1+t^{2})^{-w'}e^{-2c_0|t|^{\alpha}}\,dt,\nonumber
\end{eqnarray}
where $c_4$ depends only on $C,C',w',\pi$. Combining
(\ref{interm-minmax-1})--(\ref{minmax-interm3}) and the explicit
formula for the support of the Meyer wavelet, we can bound
$K(P^{n}_{k}|P^{n}_{0})$ by
\[
c_{5}n\gamma_{j}^{2}2^{-j}\biggl(\int_{(2\pi/3)2^{j}}^{(8\pi/3)2^{j}}(1+t^{2})^{-w}e^{-2c_0|t|^{\alpha}}\,dt
+
\int_{(2\pi/3)2^{j}}^{(8\pi/3)2^{j}}(1+t^{2})^{-w'}e^{-2c_0|t|^{\alpha}}\,dt\biggr),
\]
where $c_{5}>0$ depends only on $C,C',w',\pi,\|\psi\|_{1},\int_{\R}|x\psi(x)|\,dx$. It remains to estimate the size of these integrals and
select $j$ appropriately and we distinguish the moderately and severely
ill-posed cases.

In the moderately ill-posed case ($c_0 = 0$, $w'\ge w\ge 0$), we have
$K(P^{n}_{k}|P^{n}_{0})\le c_{6}(c')^{2}n2^{-j(2s+2w+1)}$
for some constant $c_{6}>0$ independent of $n,j$. Taking $2^{j}\simeq(n/\log n)^{1/(2s+2w+1)}$ and $c'>0$ small enough (independent of $n$
and $j$) in the definition of $\gamma_{j}$ gives
$K(P^n_{k}|P^n_{0})\le c_{6}(c')^{2}(\log n)\le\frac{1}{16}\log(\operatorname{card}(\mathcal{C}))$, where we recall that
$\operatorname{card}(\mathcal{C})=2^j$. The separation rate $r_n$ for this
choice of $j_n$ becomes, for any $k,k'$ distinct,
\[
\|f_{k}-f_{k'}\|_{\infty}\ge c_7\biggl(\frac{\log n}{n}\biggr)^{s/(2s+2w+1)}:=r_n
\]
for some constant $c_7>0$ independent of $n$. This proves Theorem
\ref{theo-minmax} for the moderately ill-posed case.

For the severely ill-posed case ($c_0>0$), we similarly obtain that
$K(P^{n}_{k}|P^{n}_{0})\le c_{8}(c')^{2}n2^{-jc(s,w,w')}2^{-d_0 2^{j\alpha}}$ with $d_0=(2 c_0 (2\pi/3)^\alpha)/\log 2$ and
constants $c_8>0$, $c(s,w,w')$ independent of $n,j$. Taking $j_n
\alpha=\log_{2}(\frac{\nu}{d_0}\log_{2}n)$ with $\nu> 1$ large enough
gives
\[
K(P^{n}_{k}|P^{n}_{0})\le c_9 (c')^2(\log_2 n)^{c'(s,w,w')}n^{1-\nu}\le\tfrac{1}{16}\log(\operatorname{card}(\mathcal{C})),
\]
where $c_9>0$, $c'(s,w,\alpha)$ are  nonnegative constants independent
of $j,n$. For this choice of $j_n,$ the separation rates $r_n$ become,
for any $k,k'$ distinct,
\[
\|f_{k}-f_{k'}\|_{\infty}\ge c_{10}\biggl(\frac{1}{\log n}\biggr)^{s/\alpha},
\]
where $c_{10}>0$ is independent of $n$. This concludes the proof of the
theorem.

\subsection{Proofs of VC properties}
\mbox{}
\begin{pf*}{Proof of Lemma \ref{entropy}}
Set
\[
\eta_j(x)= F^{-1}\biggl(1_{[-2^ja,2^ja]}\frac{1}{\overline{F[\varphi]}}\biggr)(x),
\]
which is bounded and continuous, and rewrite
\begin{eqnarray*}
\tilde\phi_{jk}(x)
&=&
\phi_{0k}(2^j\cdot)\ast\eta_j (x)
\\
&=&
\int_\mathbb R\phi(2^jx-2^jy-k)\eta_j(y)\,dy
\\
&=&
\int_\mathbb R 2^{-j/2}\phi_{j0}(x-y-2^{-j}k)\eta_j(y)\,dy
\\
&=&
2^{-j/2}\phi_{j0}\ast\eta_j(x-2^{-j}k)
\end{eqnarray*}
so that it is sufficient to study the class consisting of translates of
the fixed function $2^{-j/2}\phi_{j0}\ast\eta_j$. First, note that
$\delta_j\tilde\phi_{jk}$, $k\in\mathbb Z$, is uniformly bounded in
view of the last estimate and since
\begin{equation} \label{enve}
(2^{-j/2}\delta_j)\|\phi_{j0}\ast\eta_j\|_\infty
\le
(2^{-j/2}\delta_j)\|\phi_{j0}\|_2\|\eta_j\|_2
\le
\sqrt{2a}/2\pi,
\end{equation}
where we have used Young's convolution inequality and Plancherel's
theorem.

To prove the entropy bound, we will show that $\phi_{j0}\ast\eta_j$
has finite quadratic variation (i.e., $2$-variation). In fact, to
obtain a bound on the quadratic variation that is independent of $j$,
we renormalize and show that the function $(2^{-j/2}\delta_j)\phi_{j0}\ast\eta_j$ has quadratic variation bounded by a constant $D$ that
depends only on $\phi$. This will complete the proof of the lemma by
using Lemma 1 in \cite{GN09a}, which states that the set of dilates and
translates of a fixed function $h$ of bounded $p$-variation, $1\le p<\infty$, is of VC-type with constants $A,v$ depending only on $p$ and
the $p$-variation norm of $h$.

We will prove that $(2^{-j/2}\delta_j)\phi_{j0}\ast\eta_j$ has
bounded quadratic variation by showing that it is contained in the
(homogeneous) Besov space $\dot{B}^{1/2}_{21}(\mathbb R)$, which is
sufficient, in view of the continuous embedding of
$\dot{B}^{1/2}_{21}(\mathbb R)$ into the space $V^2(\mathbb R)$ of
functions of quadratic variation (a result due to Peetre---see Theorem
5 in \cite{BCS06} for a proof,~also the proof of Theorem 2 in
\cite{NicklPotscher2007}, which applies to $p=2$ as well). The seminorm
$\|\cdot\dot{\|}_{1/2,2,1}$ of $\dot{B}^{1/2}_{21}(\mathbb R)$ has the
following Littlewood--Paley characterization:
\[
\|h\dot{\|}_{1/2,2,1}=\sum_{l\in\mathbb Z}2^{l/2}\|F^{-1}[\gamma_l F[h]]\|_2,
\]
where $\gamma_l$ is a dyadic partition of unity with $\gamma_l$ supported in $[2^{l-1},2^{l+1}]$ (see,
e.g., Theorem 6.3.1 and Lemma 6.1.7 in \cite{BL76}). We bound the Littlewood--Paley norm:
using\vspace*{1.5pt}  the fact that $F[2^{-j/2}\phi_{j0}] = 2^{-j} F[\phi](2^{-j}\cdot)$ and Plancherel's
theorem,  introducing the notation $\langle u\rangle = (1+|u|^2)^{1/2}$ and in view of the
support of $\gamma_l$, we have the bound
\begin{eqnarray*}
&&
\delta_j\sum_l 2^{l/2}\bigl\|F^{-1}\bigl[\gamma_l F[2^{-j/2}\phi_{j0}\ast\eta_j]\bigr]\bigr\|_2
\\
&&\qquad=
\frac{1}{2\pi} 2^{-j}\delta_j\sum_l 2^{l/2}\biggl\|\gamma_l F[\phi](2^{-j}\cdot)1_{[-2^ja,2^ja]}(\overline{F[\varphi]})^{-1}\frac{\langle u\rangle^{1/2}}{\langle u\rangle^{1/2}}\biggr\|_2
\\
&&\qquad\le
c2^{-j}\delta_j\sum_l\bigl\|\gamma_l F[\phi](2^{-j}\cdot)1_{[-2^ja,2^ja]}(\overline{F[\varphi]})^{-1}\langle u\rangle^{1/2}\bigr\|_2
\\
&&\qquad\le
c 2^{-j}\sum_l\sqrt{\int_{-2^ja}^{2^ja}\gamma^2_l(u)|F[\phi](2^{-j}u)|^2\langle u\rangle\,du}
\\
&&\qquad\le
c(a)2^{-j/2}\sum_l\|F^{-1}[\gamma_l F[\phi](2^{-j}\cdot)]\|_2
\\
&&\qquad=
c(a)\sum_l\|F^{-1}[\gamma_l F[\phi_{j0}]]\|_2\le c(a)\|\phi_{j0}\dot{\|}_{0,2,1}.
\end{eqnarray*}
To bound the last quantity, we use the inequality
$\|\cdot\dot{\|}_{0,2,1}\le\|\cdot\|_{0,2,1}$  (which follows from
Definition \ref{besov} and results in \cite{M92}, Section 6.10). By
orthogonality of the wavelet basis ($j\in\mathbb N$, without loss of
generality),
\[
\|\phi_{j0}\|_{0,2,1}
=
\sqrt{\sum_k|\langle\phi_{j0},\phi_{0k}\rangle|^2}+\sum_{l=0}^{j-1}\sqrt{\sum_k|\langle\psi_{lk},\phi_{j0}\rangle|^2}.
\]
The first term on the right-hand side is bounded by $
\|K_{0}(\phi_{j0})\|_2\le\|\phi_{j0}\|_2\le 1$ since $K_0$ is an
$L^2$-projection. For the second term, we note, writing $\psi_k$ for
$\psi(\cdot-k)$ and using the change of variables $2^jx=u$ and
Condition \ref{wav}, that
\begin{eqnarray*}
\sum_k|\langle\psi_{lk},\phi_{j0}\rangle|^2
&=&
\sum_k\biggl(2^{l/2}2^{j/2}\int\psi_k(2^lx)\phi(2^jx)\,dx\biggr)^2
\\
&=&
\sum_k\biggl(2^{l/2}2^{-j/2}\int\psi_k(2^{l-j}u)\phi(u)\,du\biggr)^2
\\
&\le&
2^l2^{-j}\sup_k\biggl|\int\psi_k(2^{l-j}u)\phi(u)\,du\biggr|c(\psi)\|\phi\|_1
\\
&\le&
C^2(\psi,\phi)2^{l-j}
\end{eqnarray*}
so that
\[
\sum_{l=0}^{j-1}\sqrt{\sum_k|\langle\psi_{lk},\phi_{j0}\rangle|^2}\le C(\psi,\phi)2^{-j/2}\sum_{l=0}^{j-1}2^{l/2}\le C(\psi,\phi).
\]
This shows that $2^{-j/2}\delta_j\|\phi_{j0}\ast\eta_j\dot{\|}_{1/2,2,1}$ is bounded by a fixed constant that
depends only on $\phi,\psi$, which completes the proof of the entropy
bound. The proof of Lemma \ref{entropy1} is the same (in fact, it is
simpler since, in the last step, by orthogonality, only the resolution
level $l$ has to be considered).
\end{pf*}

\subsection{\texorpdfstring{Proofs of Propositions \protect\ref{suprat} and
\protect\ref{suprat-beta}}{Proofs of Propositions 1 and 2}}
\mbox{}
\begin{pf*}{Proof of Proposition \ref{suprat}}
We recall (\ref{sup}) and observe that $\mathcal H_j$ is bounded by the
fixed constant $U$. We prove $j>0$; the case $j=0$ is the same, except
for notation. Using the moment inequality (57) in \cite{GN09a} and
Lemma \ref{entropy}, we obtain
\begin{eqnarray*}
E\sup_{k\in\mathbb Z}\Biggl|\frac{2^j}{n}\sum_{m=1}^n\bigl(\tilde \phi_{jk}(Y_m)-E\tilde \phi_{jk}(Y)\bigr)\Biggr|
&=&
\frac{2^j}{\delta_j n}E\Biggl\|\sum_{m=1}^n\bigl(h(Y_m)-Eh(Y)\bigr)\Biggr\|_{\mathcal H_j}
\\
&\le&
\frac{C(v)2^{j}}{\delta_j n}\biggl(\sigma\sqrt{n\log\frac{AU}{\sigma}}+\log\frac{AU}{\sigma}\biggr)
\\
&\le&
\frac{C(v,A,U)}{\delta_j}\biggl(\sqrt{G^2\frac{2^jj}{n}}+\frac{2^j j}{n}\biggr),
\end{eqnarray*}
where $\sigma^2\ge\sup_{h\in\mathcal H_j}Eh^2(Y)$ is obtained as
follows: using Plancherel's theorem,
\begin{eqnarray*}
Eh^2(Y)
&=&
\delta_j^2\int_\mathbb R\tilde\phi_{jk}^2(x)g(x)\,dx\le\delta_j^2\|g\|_\infty\|\tilde\phi_{jk}\|_2^2
\\
&=&
\frac{1}{2\pi}\delta_j^2 2^{-2j}\|g\|_\infty\int_{-2^ja}^{2^ja}|F[\phi_{0k}](2^{-j}u)|^2|F[\varphi](u)|^{-2}\,du
\\
&\le&
\frac{1}{2\pi}2^{-2j}\|g\|_\infty\int_{-2^ja}^{2^ja}|F[\phi_{0k}](2^{-j}u)|^2\,du
\\
&\le&
\frac{1}{2\pi}2^{-j}\|g\|_\infty\int_{-a}^a|F[\phi_{0k}](v)|^2\,dv
\\
&=&
2^{-j}\|g\|_\infty\le 2^{-j}G^2\equiv\sigma^2,
\end{eqnarray*}
a bound which does not depend on $h$. The claim for general $p$ follows
from standard arguments for uniformly bounded empirical processes,
using, for instance, Proposition 3.1 in \cite{GLZ00}.

We now prove the second statement. For every $u'>0$, Talagrand's
inequality in Bousquet's version \cite{bousquet2003} applied to
$Z=\frac{2^{j}}{\delta_{j}n}\|\sum_{m=1}^n(h(Y_m)-Eh(Y))\|_{\mathcal H_j}$ yields
\[
\Pr\biggl\{Z\ge EZ+\sqrt{\frac{2u'}{\delta_{j}}\biggl(G^2\frac{2^{j}}{n\delta_{j}}+\frac{2U2^{j}}{n}EZ\biggr)}+\frac{U2^{j}u'}{3\delta_{j}n}\biggr\}\le e^{-u'}.
\]
Now, the first statement of the proposition and taking $u'=(1+u)j'$
imply, after some elementary computations, that
\begin{eqnarray*}
&&
\Pr\biggl\{\sup_{x\in\mathbb R}|f_n(x,j)-Ef_n(x,j)|\ge\frac{C}{\delta_{j}}\biggl(G\sqrt{\frac{2^jj'(1+u)}{n}}+\frac{2^j j'(1+u)}{n}\biggr)\biggr\}
\\
&&\qquad\le
e^{-(1+u)j'},
\end{eqnarray*}
which completes the proof.
\end{pf*}

\begin{pf*}{Proof of Proposition \ref{suprat-beta}}
The proof is the same as that of Proposition \ref{suprat} (up to some
obvious modifications).
\end{pf*}

\subsection{\texorpdfstring{Proofs of Theorems \protect\ref{uni-band} and
\protect\ref{mod-ill-theo}}{Proofs of Theorems 2 and 3}}

First, consider Theorem \ref{uni-band}. The bias is
\[
\sup_{x\in\mathbb R}|f(x)-Ef_{n}(x,j_{n})|
=
\|f-K_{j_n}(f)\|_\infty
\le
C_{1}2^{-j_{n}s}
\le
C_{1}'\biggl(\frac{1}{\nu\log n}\biggr)^{s/\alpha},
\]
where $C_{1}'>0$ depends only on $\|f\|_{s,\infty,\infty}$ (see
Theorem 9.4 in \cite{HKPT}). For the ``variance'' term, Proposition
\ref{suprat} and our choice for $j_{n}$ give
\begin{eqnarray*}
&&E\sup_{x\in\mathbb R}|f_{n}(x,j_{n})-Ef_{n}(x,j_{n})|
\\
&&\qquad\le
L''''e^{c_0 a^\alpha 2^{j_{n}\alpha}}\biggl(G \sqrt{(\nu\log n)^{1/\alpha}\frac{\log_{2}(\nu\log n)}{\alpha n}}+(\nu\log n)^{1/\alpha}\frac{\log_{2}(\nu\log n)}{\alpha n}\biggr)
\\
&&\qquad\le
L'''''G\frac{n^{c_0 a^\alpha\nu}}{\sqrt n}\sqrt{(\log n)^{1/\alpha}\log\log n}=o\biggl(\biggl(\frac{1}{\log n}\biggr)^{s/\alpha}\biggr).
\end{eqnarray*}

Using Proposition \ref{suprat} and the above bias-variance
decomposition, the proof of Theorem \ref{mod-ill-theo} is similar to
that of Theorem \ref{uni-band} and is left to the reader.

\subsection{\texorpdfstring{Proof of Theorem
\protect\ref{result-threshold}}{Proof of Theorem 4}}

For simplicity of notation, we suppress the suprema over $B(s,L)$ in
most of what follows---uniformity of the bound follows from tracking
all of the constants involved and noting that any density in $B(s,L)$
is bounded by a fixed constant $U$ that depends only on $s,L$. We have
\begin{eqnarray*}
\sup_{f\in B(s,L)}E\|f_n^{T}-f\|_{\infty}
&\leq&
\sup_{f\in B(s,L)}E\sup_{y\in\mathbb R}|f_n(y,0)-Ef_n(y,0)|
\\
&&{}+
\sup_{f\in B(s,L)}E\Biggl\|\sum_{l=0}^{j_{1}-1}\sum_{k}\bigl(\hat{\beta}_{lk}1_{|\hat{\beta}_{lk}|>\tau(l)}-\beta_{lk}(f)\bigr)\psi_{lk}\Biggr\|_{\infty}
\\
&&{}+
\sup_{f\in B(s,L)}\|K_{j_{1}}(f)-f\|_{\infty}.
\end{eqnarray*}
The first term in the right-hand side is treated in Proposition \ref{suprat},
which implies that $\sup_{f\in B(s,L)} E \sup_{y \in \mathbb R}
|f_n(y,0)-Ef_n(y,0)| \leq c\sqrt{1/n},$ which is of smaller
order than the right-hand side~in (\ref{maint}). For the third,
``deterministic,'' term, we have, from standard approximation results
for wavelets (Theorem 9.4 in \cite{HKPT}), $ \| K_{j_{1}}(f) - f
\|_{\infty} \leq c(L)  2^{-j_{1}s} \le c'(L)((\log n)/n
)^{s/(2w+1)},$ which is again of smaller order than the
right-hand side~in (\ref{maint}).

The quantity inside  the expectation of the supremum of the second term
can be decomposed, for any $f\in B(s,L)$, as
\begin{eqnarray*}
&&
\sum_{l=0}^{j_{1}-1}\sum_{k}(\hat{\beta}_{lk}-\beta_{lk})\psi_{lk}\bigl(1_{|\hat{\beta}_{lk}|>\tau,|\beta_{lk}|>\tau/2}+1_{|\hat{\beta}_{lk}|>\tau,|\beta_{lk}|\leq\tau/2}\bigr)
\\
&&\qquad{}-
\sum_{l=0}^{j_{1}-1}\sum_{k}\beta_{lk}\psi_{lk}\bigl(1_{|\hat{\beta}_{lk}|\le\tau,|\beta_{lk}|>2\tau}+1_{|\hat{\beta}_{lk}|\le\tau,|\beta_{lk}|\leq 2\tau}\bigr)
\end{eqnarray*}
and we denote these terms (I)--(IV).

We first treat the ``large deviation'' terms (II) and (III). For (II),
using (\ref{abs}) and the  Cauchy--Schwarz inequality, we have
\begin{eqnarray}\label{interm1}
&&
E\sup_{y\in\R}\Biggl|\sum_{l=0}^{j_{1}-1}\sum_{k}(\hat{\beta}_{lk}-\beta_{lk})\psi_{lk}(y) 1_{|\hat{\beta}_{lk}|>\tau,|\beta_{lk}|\leq\tau/2}\Biggr|\nonumber
\\
&&\qquad\leq
E\Biggl[\sum_{l=0}^{j_{1}-1}\sup_{k}|\hat{\beta}_{lk}-\beta_{lk}|\sup_{k}1_{|\hat{\beta}_{lk}|>\tau,|\beta_{lk}|\leq\tau/2}\sup_{y\in\R}\sum_{k}|\psi_{lk}(y)|\Biggr]
\\
&&\qquad\leq
\sum_{l=0}^{j_{1}-1}2^{l/2}c(\psi)\Bigl[E\sup_{k}|\hat{\beta}_{lk}-\beta_{lk}|^{2}\Bigr]^{1/2}\Bigl[E\sup_{k}1_{|\hat{\beta}_{lk}|>\tau,|\beta_{lk}|\leq\tau/2}\Bigr]^{1/2}.\nonumber
\end{eqnarray}
We have, using the second part of Proposition \ref{suprat-beta}, choosing $\kappa'$
large enough  depending only on $a,w,C,\phi,\psi$ and using the fact
that $(2^l l/n)^{1/2}$ is bounded by a fixed constant independent of
$l$,
\begin{eqnarray}\label{interm2}
&&
E\sup_{k}1_{|\hat{\beta}_{lk}|>\tau,|\beta_{lk}|\leq\tau/2}\nonumber
\\
&&\qquad\leq
E\Bigl(\sup_{k}1_{|\hat{\beta}_{lk}-\beta_{lk}|>\tau/2}\Bigr)\nonumber
\\
&&\qquad\le
\Pr\biggl(\sup_{k}|\hat{\beta}_{lk}-\beta_{lk}|>\frac{\kappa'}{2a^w} G 2^{lw}a^w \sqrt{\frac{\log n}{n}}\biggr)\nonumber
\\[-8pt]\\[-8pt]
&&\qquad\le
\Pr\biggl(\sup_{k} |\hat{\beta}_{lk}-\beta_{lk}|>c(a,w,C)\kappa' G\frac{1}{\delta_l}\sqrt{\frac{\log n}{n}}\biggr)\nonumber
\\
&&\qquad\le
\Pr\biggl(\sup_{k} |\hat{\beta}_{lk}-\beta_{lk}|>\frac{c(a,w,C)\kappa'}{\delta_l}G\sqrt{\frac{[1+(\log n/l')-1]l'}{n}}\biggr)\nonumber
\\
&&\qquad\le
e^{-2\log n}.\nonumber
\end{eqnarray}
Now, combining (\ref{interm1}) and (\ref{interm2}) with the first part
of Proposition \ref{suprat-beta} yields the bound
\begin{eqnarray*}
c\sum_{l=0}^{j_{1}-1}2^{l(w+(1/2))}G\sqrt{\frac{l'}{n}} e^{-\log n}
&\le&
C' G e^{-\log n} \sqrt{\frac{\log n}{n}} 2^{j_1(w+1/2)}
\\
&\le&
\frac{C''}{n} = o(n^{-1/2})
\end{eqnarray*}
for (II).

For term (III), using (\ref{interm2}), as well as
$\sum_{k}|\beta_{lk}|\leq c(\psi) 2^{l/2}$ for any density $f$, we have
\begin{eqnarray*}
&&
E\sup_{y\in \R}\Biggl|\sum_{l=0}^{j_{1}-1}\sum_{k}\beta_{lk}\psi_{lk}(y)1_{|\hat{\beta}_{lk}|\leq \tau, |\beta_{lk}|>2\tau}\Biggr|
\\
&&\qquad\leq
\sum_{l=0}^{j_{1}-1}2^{l/2}\|\psi\|_{\infty}\sum_{k}|\beta_{lk}|\Pr(|\hat{\beta}_{lk}|\leq\tau, |\beta_{lk}|>2\tau)
\\
&&\qquad\leq
C''' e^{-2\log n} \sum_{l=0}^{j_{1}-1}2^{l}  \le  C''''n^{-2} (n/\log n)^{1/(2w+1)}= o(n^{-1/2}).
\end{eqnarray*}
We now bound (I). Let $j_{1}(s)$ be such that $0\le j_{1}(s)\leq
j_{1}-1$ and
\begin{equation} \label{opt}
2^{j_{1}(s)}\simeq (n/\log n)^{1/(2s+2w+1)}
\end{equation}
[such $j_{1}(s)$ exists by the definitions]. Proposition \ref{suprat-beta} and
(\ref{abs}) give
\begin{eqnarray*}
&&
E\sup_{y\in\mathbb R}\Biggl|\sum_{l=0}^{j_{1}(s)-1}\sum_{k}(\hat{\beta}_{lk}-\beta_{lk})\psi_{lk}(y)1_{|\hat{\beta}_{lk}|>\tau,|\beta_{lk}|>\tau/2}\Biggr|
\\
&&\qquad\leq
\sum_{l=0}^{j_{1}(s)-1}E\sup_{k}|\hat{\beta}_{lk}-\beta_{lk}|2^{l/2}c(\psi)
\\
&&\qquad\leq
D G\sum_{l=0}^{j_{1}(s)-1}2^{lw}\sqrt{\frac{2^{l}l'}{n}}
\\
&&\qquad\leq
D'G 2^{j_1(s)w}\sqrt{\frac{2^{j_{1}(s)}j_{1}(s)}{n}}\leq D'' G\biggl(\frac{\log n}{n}\biggr)^{s/(2s+2w+1)},
\end{eqnarray*}
where $D''>0$ depends only on $\psi,\phi,C,w$. For the second part of
(I), using the fact that Definition \ref{besov} implies
\begin{equation} \label{biascof}
\sup_{k}|\beta_{lk}(f)|\leq D(L) 2^{-l(s+1/2)}
\end{equation}
for $f \in B(s,L)$, the definition of $\tau$ and Proposition \ref{suprat-beta}, we
obtain
\begin{eqnarray*}
&&
E\sup_{y\in\R}\Biggl|\sum_{l=j_{1}(s)}^{j_{1}-1}\sum_{k}(\hat{\beta}_{lk}-\beta_{lk})\psi_{lk}(y)1_{|\hat{\beta}_{lk}|>\tau,|\beta_{lk}|>\tau/2}\Biggr|
\\
&&\qquad\leq
\sum_{l=j_{1}(s)}^{j_{1}-1}E\sup_{k}|\hat{\beta}_{lk}-\beta_{lk}|\frac{2}{\kappa}2^{-lw}\sqrt{\frac{n}{\log n}}\sup_{k}|\beta_{lk}|2^{l/2}c(\psi)
\\
&&\qquad\leq
D'''\sum_{l=j_{1}(s)}^{j_{1}-1}2^{-ls}\le D''''\biggl(\frac{\log n}{n}\biggr)^{s/(2s+2w+1)},
\end{eqnarray*}
where $D''''$ depends only on $L,s,\kappa',\phi,\psi,C$.

To complete the proof, we control the term (IV). Again using
(\ref{biascof}), we have
\begin{eqnarray}\label{interm-(IV)}
&&\sup_{y\in\mathbb R}\Biggl|\sum_{l=0}^{j_{1}-1}\sum_{k}\beta_{lk}\psi_{lk}(y)1_{|\hat{\beta}_{lk}|\leq\tau,|\beta_{lk}|\leq
2\tau}\Biggl|\nonumber
\\
&&\qquad\leq
c(\psi)\sum_{l=0}^{j_{1}-1}\sup_{k}2^{l/2}|\beta_{lk}|1_{|\beta_{lk}|\leq 2\tau}
\\
&&\qquad\leq
c'\sum_{l=0}^{j_{1}-1}\min\biggl(2^{l(w+1/2)}\sqrt{\frac{\log n}{n}},2^{-ls}\biggr).\nonumber
\end{eqnarray}
Since the antagonistic terms in the minimum are strictly monotone in
$l$, the $l^* \in \mathbb R$ for which they are maximal is the one
where they are equal so that $2^{l^*} \simeq 2^{j_1(s)}$
[cf.~(\ref{opt})]. If we denote by $[l^*]$ the integer part of $l^*,$
then the last sum is bounded by
\[
c'\sum_{l=0}^{[l^*]}2^{l(w+1/2)}\sqrt{\frac{\log n}{n}}+c'\sum_{l=[l^*]+1}^{j_1-1} 2^{-ls}
\le
c''\biggl(\frac{\log n}{n}\biggr)^{s/(2w+2s+1)}.
\]

\subsection{\texorpdfstring{Proofs for Section \protect\ref{app}}{Proofs for Section 2.4}}

The following proposition is the wavelet-analog of a similar result in
Proposition 1 in \cite{BT08I} for kernel regularizations.

\begin{proposition} \label{smooth-bias}
Let $\phi,\psi$ satisfy Condition \textup{\ref{wav}}. Let $f\in \mathcal
A_{\tilde{c}_0,s}(L)$ for some $\tilde{c}_0,s,L>0$. We then have, for
every $j \ge 0,$ that
\[
\|K_{j}(f)-f\|_{\infty}\le c'''\sqrt{L}2^{j(1-s)/2}e^{-\tilde{c}_0(a')^{s} 2^{js}},
\]
where the constant $c'''>0$ depends only on $\phi,\psi,\tilde c_0,s$.
\end{proposition}

\begin{pf}
Using (\ref{abs}), Plancherel's theorem and the fact that $f\in
\mathcal{A}_{\tilde c_0,s}(L),$ we have
\begin{eqnarray*}
\|K_j(f)-f\|_\infty
&\le&
c(\psi) \sum_{l\ge j}2^{l/2}\sup_{ k \in \mathbb Z} |\beta_{lk}(f)|
\\
&=&
c' \sum_{l\ge j}2^{l/2}\sup_{k\in\mathbb Z}\biggl|\int_{\R}\overline{F[\psi_{lk}](u)}Ff(u)\,du\biggr|
\\
&\le&
c' \sum_{l\ge j}\sup_{k\in \mathbb Z} \int_{\R} |F[\psi](2^{-l}u)| |Ff(u)|\,du
\\
&\le&
c'\|\psi\|_1 \sum_{l\ge j}\int_{\mathbb{R}\setminus[-2^la',2^la']}|Ff(u)|e^{\tilde c_0 |u|^s}e^{-\tilde c_0|u|^s}\,du
\\
&\le&
c''\|\psi\|_1 \sqrt{L}\sum_{l\ge j}\sqrt{\int_{2^la'}^{\infty}e^{-2\tilde c_0 u^s}\,du}
\end{eqnarray*}
and the result follows from the inequality $\int_a^\infty e^{-cu^s}\,du\le C(c,s)a^{1-s} e^{-ca^s}$ for $a,s>0$.
\end{pf}

\begin{pf*}{Proof of Corollary \ref{supsm}}
Decomposing the sup-norm error of the linear estimator into ``bias'' and
``variance'' terms and applying Propositions \ref{suprat} and \ref{smooth-bias}, we have, for any $j\ge 0$,
\begin{eqnarray*}
&&E\sup_{x\in\mathbb{R}}|f_n(x,j)-f(x)|
\\
&&\qquad\le
\sup_{x\in\mathbb{R}}|Ef_n(x,j)-f(x)|+E\sup_{x\in\mathbb{R}}|f_n(x,j)-E f_n(x,j)|
\\
&&\qquad\le
c'''\sqrt{L}e^{-\tilde{c}_0(a')^{s}2^{js}} 2^{j(1-s)/2}+c\frac{1}{\delta_j}\biggl(G\sqrt{\frac{2^jj'}{n}}+\frac{2^j j'}{n}\biggr)
\\
&&\qquad\le
C'\biggl(e^{-\log n/2}(\log n)^{(1-s)/2s}+ 2^{jw}\sqrt{\frac{2^j j'}{n}}\biggr),
\end{eqnarray*}
where $C'>0$ depends only on $C,s,L,\tilde c_0,a,w$. The result
follows immediately for $s\ge 1$ and for $s<1$ in view of the fact
that $(1-s)/2s<(w+1/2)/s$ for all $s>0$.
\end{pf*}

\begin{pf*}{Proof of Theorem \ref{minm2}}
The proof of this theorem follows the one of Theorem \ref{theo-minmax}
up to the following modifications. Let $p$ be the standard Cauchy
density. Fix $0<\nu<1/2$. Since $F[p](u)=e^{-|u|}$, we see from the
scaling property of Fourier transforms and since $s \le 1$ that there
exists a constant $\eta=\eta(\nu)$ large enough such that
$f_0=(1/\eta)p(\cdot/\eta) \in \mathcal{A}_{\tilde{c}_0,s}(\nu^2L)$.

As in the proof of Theorem \ref{theo-minmax}, we consider the functions
$ f_{k}(x) = f_0(x) + \gamma_j \psi_{jk_M}$, $1 \le k \le 2^j -1$,
$k_M=kM$, $M \ge 1$ with $\gamma_j =
c'\sqrt{L}\sqrt{j}2^{jw}e^{-\tilde{c}_0[a^s+1]2^{js}}.$ We have $f_k\in
\mathcal{A}_{\tilde{c}_0,s}(L)$ for every $k$ if $c'>0$ is a constant
taken small enough and depending only on $\nu,a,\|\psi\|_1$ since
\begin{eqnarray*}
&&
\int_{\mathbb R}|F[f_k](t)|^2e^{2\tilde{c}_0|t|^s}\,dt
\\
&&\qquad\le
2\int_{\mathbb R}|F[f_0](t)|^2e^{2\tilde{c}_0|t|^s}dt+2\gamma_j^2\int_{\mathbb R}|F[\psi_{j,k}](t)|^2e^{2\tilde{c}_0|t|^s}\,dt
\\
&&\qquad\le
4\pi\nu^{2}L+2\gamma_j^22^{-j}\|\psi\|_1^2\int_{a'2^{j}}^{a2^{j}}e^{2\tilde{c}_0|t|^s}\,dt
\\
&&\qquad\le
4\pi\nu^{2}L+2(c')^2Lj2^{2jw}e^{-2\tilde{c}_0[a^s+1]2^{js}}\|\psi\|^2_1 a2^je^{2\tilde{c}_0a^s2^{js}}
\\
&&\qquad\le
2\pi L.
\end{eqnarray*}
Take $2^{js} = \frac{1}{2\tilde{c}_0[a^s+1]}\log n$. The proof of
Theorem \ref{theo-minmax} then implies, $\forall k\neq k'$, $ \|f_k -
f_{k'}\|_{\infty} \ge c_3 \sqrt{(\log \log n)/n}(\log n)^{(w+1/2)/s}$
for some constant $c_3>0$ independent of $n$. Next, for any $k,$ the
Kullback--Leibler divergence between $P^n_{k}$ and $P^n_{0}$ satisfies
\[
K(P^n_{k}|P^n_{0}) \le c_4 n \gamma_j^2 2^{-2jw} = c_{4}(c')^{2}L n j 2^{2jw} e^{-2\tilde{c}_0[a^s+1]2^{js}}2^{-2jw}\le c_{4}(c')^{2}L j.
\]
This and Lemma \ref{lemlb} together yield the
result for $c'>0$ chosen small enough independently of $n, k$.
\end{pf*}

\begin{pf*}{Proof of Proposition \ref{rad}}
We use Proposition \ref{rad0} below. Note that
\begin{equation} \label{sup0}
\|f_n(j)-Ef_n(j)\|_\infty
=
\sup_{x \in \mathbb R}\Biggl|\frac{1}{n}\sum_{m=1}^n\bigl(K^*_j(x,Y_m)-EK^*_j(x,Y)\bigr)\Biggr|.
\end{equation}
The class $\{K^*_j(x, \cdot)\dvtx x \in \R \}$ has\vspace*{1pt} envelope $U(j)=2^j
\delta_j^{-1}c(\phi)\sqrt{a/2\pi^2}$ in view of (\ref{abs}) and
(\ref{enve}). Since Proposition \ref{rad0} deals with classes of
functions bounded by $1/2,$ we have to rescale, that is, we consider
the class $\mathcal G := \mathcal G_j = \{K^*_j(x, \cdot)/\break 2U(j)\dvtx x \in
\R \},$ which is uniformly bounded by $1/2$. Furthermore, the upper
bound for the weak variances $\sup_{g \in \mathcal G}Eg^2(Y) \le
\sigma^2$ can be taken to be $2^{-j}(\pi/2) \|\phi\|^2_1\|g\|_\infty$
in view of the estimate
\begin{eqnarray*}
E(K^*_j(x,Y))^2
&\le&
 2^j\|g\|_\infty c(\phi)^2
\|\phi_{j0}\|_1^2 \|\eta_j\|_2^2
\\
&\le&
 \|g\|_\infty c(\phi)^2
\|\phi\|_1^2 \delta_j^{-2} 2^j (a/\pi),
\end{eqnarray*}
which uses Young's inequality
(and the definition of $\eta_j$ from the proof of Lemma~\ref{entropy}).

To prove the inequality, set $d(\phi)=c(\phi)\sqrt{a/2 \pi^2}$ and
$d'(\phi)=d(\phi)\|\phi\|_1\sqrt{2\pi}$ so that
\begin{eqnarray*}
&&\Pr\biggl\{\|f_n(j,\cdot)-Ef_n(j,\cdot)\|_\infty
\\
&&\hphantom{\Pr\biggl\{}\ge
6R_n(j)+\frac{10d'(\phi)}{\delta_j}\sqrt{\frac{2^{j}\|g\|_\infty(z+ \log 2)}{n}}+\frac{44}{\delta_j}\frac{2^{j}d(\phi)(z+\log2)}{n}\biggr\}
\\
&&\qquad=
\Pr\Biggl\{\Biggl\|\frac{1}{n}\sum_{m=1}^n\frac{(K^*_j(\cdot,Y_m)-EK^*_j(\cdot,Y))}{2U(j)}\Biggr\|_\infty
\\
&&\qquad\hphantom{=
\Pr\Biggl\{}\ge
\frac{6R_n(j)}{2U(j)} + 10\|\phi\|_1\sqrt{\frac{\pi\|g\|_\infty(z+\log 2)}{2^{j+1}n}} + 22 \frac{z+\log 2}{n}\Biggr\},
\end{eqnarray*}
but this quantity equals the probability in
Proposition \ref{rad0} below for $\mathcal F =\mathcal G$.

For the second claim of the proposition, we only have to show that
$ER_n(j)$ has, up to constants, the required order as a function of
$j,n$. But this follows readily from the usual desymmetrization
inequality for Rademacher processes (cf., e.g., expression (23) in
\cite{GN10b}), as well as from Proposition \ref{suprat}.
\end{pf*}

\begin{pf*}{Proof of Corollary \ref{band}}
The result follows from standard arguments (combining Propositions \ref{rad} and \ref{smooth-bias}).
\end{pf*}

\subsection{A concentration inequality using Rademacher processes}

We start with the following inequality, which is a Bernstein-type
version of similar inequalities in \cite{K06} and complements the results in \cite{GN10b}. Let $\|H\|_\mathcal F =
\sup_{f \in \mathcal F}|H(f)|$ for any set $\mathcal F$ and functions
$H: \mathcal F \to \mathbb R$.

\begin{proposition}\label{rad0}
Let $X_1,\ldots,X_n$ be i.i.d.~with law $P$ on a measurable
space $(S, \mathcal A)$. Let $\mathcal F$ be a countable class of real-valued measurable functions
defined on $S$, uniformly bounded by $1/2$, and let $\sigma^2\ge\sup_{f\in\mathcal F}Ef^2(X)$. We have, for
every $n\in\mathbb N$ and $x>0$, that $e^{-x}$ is greater than or equal to
\begin{eqnarray*}
&&
\Pr\Biggl\{\Biggl\|\frac{1}{n}\sum_{i=1}^n\bigl(f(X_i)-Pf\bigr)\Biggr\|_\mathcal F
\\
&&\hphantom{\Pr\Biggl\{}\ge
6\Biggl\|\frac{1}{n}\sum_{i=1}^n\varepsilon_i f(X_i)\Biggr\|_\mathcal F+10\sqrt{\frac{(x+\log 2)\sigma^2}{n}}+22\frac{x+\log 2}{n}\Biggr\}.
\end{eqnarray*}
\end{proposition}

\begin{pf}
We first recall the lower-deviation version of Talagrand's inequality,
as given in \cite{KR05}, and a simple consequence of it. Using the
notation $Z = \|\sum (f(X_{i}) - Pf)\|_{\mathcal{F}}$, we have, using
the inequalities $\sqrt{a+b} \le \sqrt a + \sqrt b$ and $\sqrt{ab} \le
(a+b)/2$, that
\begin{eqnarray*}
e^{-x}
&\ge&
\Pr\bigl\{Z\le EZ-\sqrt{2x(n\sigma^2+2 EZ)}-x\bigr\}
\\
&\ge&
\Pr\bigl\{Z\le 0.5E Z-\sqrt{2xn\sigma^2}- 3x\bigr\}
\\
&=&
\Pr\Biggl\{\Biggl\|\frac{1}{n}\sum_{i=1}^n\bigl(f(X_i)-Pf\bigr)\Biggr\|_\mathcal F
\\
&&\hphantom{\Pr\Biggl\{}\le
0.5 E\Biggl\|\frac{1}{n}\sum_{i=1}^n\bigl(f(X_i)-Pf\bigr)\Biggr\|_\mathcal F -\sqrt{\frac{2x\sigma^2}{n}}-\frac{3x}{n}\Biggr\}
\end{eqnarray*}
and  one likewise proves, using the upper-deviation version of
Talagrand's inequality \cite{bousquet2003},
\begin{eqnarray} \label{bq}
&&
e^{-x}\ge\Pr\Biggl\{\Biggl\|\frac{1}{n}\sum_{i=1}^n\bigl(f(X_i)-Pf\bigr)\Biggr\|_\mathcal F\nonumber
\\[-8pt]\\[-8pt]
&&\hphantom{e^{-x}\ge\Pr\Biggl\{}\ge
1.5E\Biggl\|\frac{1}{n}\sum_{i=1}^n\bigl(f(X_i)-Pf\bigr)\Biggr\|_\mathcal F+\sqrt{\frac{2x\sigma^2}{n}}+\frac{7x}{3n}\Biggr\}.\nonumber
\end{eqnarray}
To prove the proposition, observe that
\begin{eqnarray*}
&&
\Pr\biggl\{\biggl\|\frac{1}{n}\sum \bigl(f(X_i)-Pf\bigr)\biggr\|_\mathcal F
\ge
6\biggl\|\frac{1}{n}\sum \varepsilon_i f(X_i)\biggr\|_\mathcal F + 10\sqrt{\frac{x\sigma^2}{n}} + \frac{22 x}{n} \biggr\}
\\
&&\qquad\le
\Pr\biggl\{\biggl\|\frac{1}{n}\sum \bigl(f(X_i)-Pf\bigr)\biggr\|_\mathcal F
\\
&&\qquad\hphantom{\le\Pr\biggl\{}\ge
3E\biggl\|\frac{1}{n}\sum \varepsilon_i f(X_i)\biggr\|_\mathcal F+1.5\sqrt{\frac{x\sigma^2}{n}} + 0.15\frac{22x}{n}\biggr\}
\\
&&\qquad\quad{}+
\Pr \biggl\{6\biggl\|\frac{1}{n}\sum \varepsilon_i f(X_i)\biggr\|_\mathcal F - 3E\biggl\|\frac{1}{n}\sum \varepsilon_i f(X_i)\biggr\|_\mathcal F
\\
&&\hphantom{\qquad\quad+\Pr \biggl\{6\biggl\|\frac{1}{n}\sum \varepsilon_i f(X_i)\biggr\|}<
-8.5\sqrt{\frac{x\sigma^2}{n}} - 0.85 \frac{22x}{n}\biggr\}
\\
&&\qquad\le
\Pr\biggl\{\biggl\|\frac{1}{n}\sum \bigl(f(X_i)-Pf\bigr)\biggr\|_\mathcal F
\\
&&\qquad\hphantom{\le\Pr\biggl\{}\ge
1.5E\biggl\|\frac{1}{n}\sum \bigl(f(X_i)-Pf\bigr)\biggr\|_\mathcal F +  \sqrt{\frac{2x\sigma^2}{n}} + \frac{7x}{3n}  \biggr\}
\\
&&\qquad\quad{}+
\Pr \biggl\{\biggl\|\frac{1}{n}\sum \varepsilon_i f(X_i)\biggr\|_\mathcal F
<
0.5 E\biggl\|\frac{1}{n}\sum \varepsilon_i f(X_i)\biggr\|_\mathcal F - \sqrt{\frac{2x\sigma^2}{n}} - \frac{3x}{n}\biggr\},
\end{eqnarray*}
where we have used the standard Rademacher symmetrization inequality
(e.g., (23) in \cite{GN10b}). The first quantity on the right-hand
side~of the last inequality is less than or equal to $e^{-x}$, by
(\ref{bq}). For the second term, note that the first displayed
inequality in this proof also applies to the randomized sums
$\sum_{i=1}^n \varepsilon_i f(X_i)$, by taking $\mathcal
G=\{g(\tau,x)=\tau f(x)\dvtx f\in \mathcal F \}$, $\tau \in
\{-1,1\}$, instead\vspace*{1.5pt} of $\mathcal F$ and the probability measure $\bar P
= 2^{-1}(\delta_{-1}+\delta_1)\times P$ instead of $P$. It is easy to
see that $\sigma$ can be taken to be the same as for $\mathcal F$. This
gives the overall bound $2e^{-x}$ and a change of variables in $x$
gives the final bound.
\end{pf}

\section*{Acknowledgments}
Both authors would like to thank two anonymous referees for their
remarks, and the second author gratefully acknowledges the hospitality
of the Cafe Braeunerhof in Vienna.

\printaddresses

\end{document}